\documentclass[11pt,oneside,leqno]{amsart}

\makeatletter
\@namedef{subjclassname@2020}{\textup{2020} Mathematics Subject Classification}
\makeatother

\usepackage[top=2.5 cm,bottom=2 cm,left=2 cm,right=3 cm]{geometry}
\usepackage[utf8]{inputenc}
\usepackage{color}
\usepackage{amssymb}

\usepackage[english]{babel}

\usepackage{todonotes}
\usepackage[backref=page]{hyperref}

\usepackage{esint}
\usepackage{graphicx}
\usepackage{hyperref}
\usepackage{bm}
\usepackage{xcolor}

\newcommand*{\mailto}[1]{\href{mailto:#1}{\nolinkurl{#1}}}

\numberwithin{equation}{section}

\newtheorem{example}{Example}[section]

\newtheorem{theorem}[example]{Theorem}
 \newtheorem{proposition}[example]{Proposition}
\newtheorem{lemma}[example]{Lemma}
 \newtheorem{corollary}[example]{Corollary}
\newtheorem{remark}[example]{Remark}
\newtheorem*{maintheorem*}{Main Theorem}
\allowdisplaybreaks
\numberwithin{equation}{section}

\renewcommand{\i}{\ifmmode\mathit{\mathchar"7010 }\else\char"10 \fi}
\renewcommand{\j}{\ifmmode\mathit{\mathchar"7011 }\else\char"11 \fi}

\newcommand{\ve}{\varepsilon}

{%

\begin{enumerate}}%
{\end{enumerate}}

{%

\begin{enumerate}}%
{\end{enumerate}}

\begin{document}

\title[Nonlocal FORQ equation and two-component peakon system]
{On the well-posedness of a nonlocal (two-place) FORQ equation  
via a two-component peakon system}

 \author[Karlsen]{K. H. Karlsen}
\author[Rybalko]{Ya. Rybalko}

\address[Kenneth Hvistendahl Karlsen]{\newline
  Department of Mathematics, \newline University of Oslo, \newline
  PO Box 1053, Blindern -- 0316 Oslo, Norway}
\email[]{kennethk@math.uio.no}

\address[Yan Rybalko]{\newline 
	Mathematical Division, 
	\newline B.Verkin Institute for Low Temperature Physics and Engineering
	of the National Academy of Sciences of Ukraine,
	\newline 47 Nauky Ave., Kharkiv, 61103, Ukraine}
\email[]{rybalkoyan@gmail.com}



\thanks{Yan Rybalko's research was partly funded through 
the Universities 4 Ukraine Fellowship, an award conferred 
by the University of California, Berkeley.}

\date{\today}

\begin{abstract}
We investigate the Cauchy problem for a
nonlocal (two-place) FORQ equation. By interpreting this equation 
as a special case of a two-component peakon system 
(exhibiting a cubic nonlinearity), 
we convert the Cauchy problem into a system of ordinary 
differential equations in a Banach space. 
Using this approach, we are able to demonstrate local 
well-posedness in the Sobolev space $H^{s}$ where $s > 5/2$. 
We also establish the continuity properties for the data-to-solution map 
for a range of Sobolev spaces. Finally, we briefly explore 
the relationship between the two-component system 
and the bi-Hamiltonian AKNS hierarchy.
\end{abstract}

\maketitle

\tableofcontents

\section{Introduction}
We are interested in the following nonlocal counterpart of the 
the Fokas-Olver-Rosenau-Qiao (FORQ) equation:
\begin{equation}\label{nFORQ-eq}
\begin{split}
&\partial_t m(x,t)=
\partial_x[m(x,t)(u(x,t)-\partial_xu(x,t))
(u(-x,-t)+\partial_x(u(-x,-t)))],\\
&m(x,t)=u(x,t)-\partial_{x}^2u(x,t).
\end{split}
\end{equation}
This equation, noted for its unique characteristics and 
dynamic behavior, was initially introduced within the context 
of two-place (Alice-Bob) integrable systems, as detailed in \cite{LH17}. 
The specific reference to this equation can be found 
in \cite[equation (26)]{LQ17}. We will delve into the 
background details and further contextual information 
later in the introduction.

Under the nonlocal transformation $n(x,t)=m(-x,-t)$, 
the equation \eqref{nFORQ-eq} can be viewed as a special 
case of the following two-component peakon system, initially 
introduced in \cite{SQQ11}:
\begin{subequations}\label{two-comp}
	\begin{align}
	\label{two-comp-a}
	&\partial_tm=\partial_x[m(u-\partial_xu)
	(v+\partial_xv)],
	&& m=m(x,t),\,u=u(x,t),\,v=v(x,t),\\
	\label{two-comp-b}
	&\partial_tn=\partial_x[n(u-\partial_xu)
	(v+\partial_xv)],
	&& n=n(x,t),\\
	&m=u-\partial_{x}^2u,\,\,n=v-\partial_{x}^2v,
	&& u,v\in\mathbb{R},\quad
	x\in A,\,\,t\in\mathbb{R}.
	\end{align}
\end{subequations}
In this paper, we investigate the Cauchy problem for 
this system, where $A$ is either 
the real line $\mathbb{R}$ or the circle 
$\mathbb{R}/2\pi\mathbb{Z}$, with
initial data $u(x,0)=u_0(x)$ and $v(x,0)=v_0(x)$ that 
belong to the Sobolev space $H^s(A)$, for $s>\frac{5}{2}$.
If we take $v(x,t)=u(x,t)$ in \eqref{two-comp} we obtain 
the Fokas-Olver-Rosenau-Qiao (FORQ) equation
\begin{equation*}
\partial_tm=\partial_x[m(u^2-(\partial_xu)^2)],
\quad m=u-\partial_{x}^2u,\quad u=u(x,t).
\end{equation*}
This equation, commonly referred to as the 
FORQ equation, finds its origin in fluid dynamics \cite{Q06}, as well 
as in the broader fields of integrable systems and soliton 
theory \cite{F95, OR96}. Over the years, numerous researchers have 
dedicated their efforts to extensively studying the FORQ 
equation \cite{BKS20, FGLQ13, GLL18, HM14, HFQ17, LLQ14}, 
recognizing its significance in multiple domains of physics. Applications 
of this equation span across diverse areas, including 
fluid mechanics and nonlinear optics \cite{GLOQ13, SW04}.

The nonlocal counterpart \eqref{nFORQ-eq} of the FORQ equation 
involves the consideration of solution values at non-neighboring 
points, such as $x$ and $-x$. The study of these systems has 
garnered interest due to their ability to describe phenomena 
characterized by strong correlations and entanglement between 
events occurring at different locations. Such systems provide 
valuable insights into complex physical phenomena 
that exhibit nonlocal behavior. Besides, in Section \ref{Hamilt-dual} 
of our paper, we establish a connection between \eqref{nFORQ-eq}, the 
nonlocal counterpart of the FORQ equation, and a nonlocal 
version of the modified Korteweg-de Vries (mKdV) equation. 
This relationship sheds some light on the underlying dynamics 
and of the interplay between these two nonlocal equations. 
The nonlocal version of the modified 
Korteweg-de Vries (mKdV) equation takes the form
\begin{equation}\label{nmkdv}
\partial_t\tilde m(x,t)
+\partial_x^3\tilde m(x,t)
+6\tilde m(x,t)\tilde m(-x,-t)
\partial_x\tilde m(x,t)=0.
\end{equation}
This equation, introduced in \cite{AM17}, holds  
relevance in the study of atmospheric and oceanic 
dynamical systems \cite{TLH18}. Finally, notice that the nonlocal 
FORQ equation \eqref{nFORQ-eq} admits peakon traveling 
 wave solutions \cite{LQ17} of the form
$u_c(x,t)=c\sqrt{\frac{3}{2}}e^{-|x+c^2t|}$.

Another exciting class of nonlocal integrable equations traces its 
origins back to the work of Ablowitz and Musslimani \cite{AM13}. 
In their paper, they introduced the nonlocal 
nonlinear Schrödinger (NNLS) equation:
\begin{equation}\label{NNLS}
\mathrm{i}\partial_{t}q(x,t)+\partial_{x}^2q(x,t)
+2\sigma q^{2}(x,t)\bar{q}(-x,t)=0.
\end{equation}
Here $\mathrm{i}^2=-1$, $q\in\mathbb{C}$, $\sigma=\pm1$ 
and $\bar{q}$ is a complex conjugate of $q$.
Equation \eqref{NNLS} can be derived as a 
nonlocal reduction of the member of the
Ablowitz-Kaup-Newell-Segur (AKNS) hierarchy \cite{AKNS74} 
and is viewed as a nonlocal counterpart of the 
conventional nonlinear Schr\"odinger equation
$$
\mathrm{i}\partial_tq(x,t)+\partial_{x}^2q(x,t)
+2\sigma |q(x,t)|^{2}q(x,t)=0.
$$
It is worth noting that both the NNLS equation \eqref{NNLS} and the 
nonlocal FORQ equation \eqref{nFORQ-eq} satisfy the parity-time-symmetric 
(PT-symmetric) condition \cite{BB98, BBCF07}. This condition asks 
that the system remains invariant under the combined operations 
of parity (P), time reversal (T), and complex conjugation.
In the case of the FORQ equation and the NNLS equation, if $u(x,t)$ 
and $q(x,t)$ represent solutions of \eqref{nFORQ-eq} and \eqref{NNLS}, 
respectively, then their respective counterparts $u(-x,-t)$ and $\bar{q}(-x,-t)$ 
are also solutions. PT-symmetric systems have attracted interest in modern 
physics due to their intriguing features and their occurrence 
in a wide range of physical systems, including optics, electronics, acoustics, and 
quantum mechanics (see, e.g., \cite{EMK18, KYZ16} and references therein). 
These systems have proven to be useful tools for investigating the properties 
of non-Hermitian systems, such as the scattering of light in 
non-uniform media, the behavior of open quantum systems, and the 
stability analysis of chaotic systems.  

The Cauchy problem for the two-component peakon system \eqref{two-comp} 
has been the subject of previous investigations.  In \cite{MM13}, the 
authors established the local existence and uniqueness 
of the solution within a range of Besov spaces, 
employing Danchin's arguments of approximate solutions \cite{D01}. 
Additionally, they obtained analytical properties of the 
solution by leveraging the Cauchy-Kowalevskaya theorem. 
In \cite{YQZ15} (also referenced as \cite{WY23}), the authors 
derived blow-up criteria for solutions of \eqref{two-comp}.

In this paper, we present a rigorous analysis of the Cauchy problem 
for \eqref{two-comp} in the Sobolev space $H^s(A)$, where $s>\frac{5}{2}$. 
Our main contributions are the establishment of local existence 
and uniqueness of the solution, along with demonstrating its continuous 
dependence on the initial data. Specifically, our findings can be summarized as follows:

\begin{theorem}[Local existence and uniqueness]
	\label{ex-un-sol}
	Consider the Cauchy problem for system \eqref{two-comp} 
	with initial data $u_0(x)=u(x,0)$ and $v_0(x)=v(x,0)$, $x\in A$.
	Assume that $u_0,v_0\in H^s$, $s>\frac{5}{2}$.
	Then there exists a unique local solution 
	$u,v\in
	C([-T_{\delta_0},T_{\delta_0}], H^{s})
	\cap C^1([-T_{\delta_0},T_{\delta_0}], H^{s-1})$, where 
	$T_{\delta_0}=
	\frac{1-\delta_0}
	{8C_s(\|u_0\|_{H^s}+\|v_0\|_{H^s})^2}$ with some $C_s>0$ 
	and $0<\delta_0<1$, which satisfy the following size estimate for the solution
	\begin{equation}\label{size-est-u-v}
	\|u(t)\|_{H^s}+\|v(t)\|_{H^{s}}
	\leq
	\frac{2}{\sqrt{\delta_0}}
	(\|u_0\|_{H^s}+\|v_0\|_{H^s}),\quad 
	-T_{\delta_0}\leq t\leq T_{\delta_0},
	\end{equation}
	as well as the size estimate for its derivative 
	\begin{equation}\label{deriv-est-u-v}
	\|\partial_tu(t)\|_{H^{s-1}}
	+\|\partial_tv(t)\|_{H^{s-1}}
	\lesssim_{s,\delta_0}
	(\|u_0\|_{H^s}+\|v_0\|_{H^s}),\quad 
	-T_{\delta_0}\leq t\leq T_{\delta_0}.
	\end{equation}
\end{theorem}

\begin{theorem}
	[Continuity of the data-to-solution map]
	\label{cont-dep}
	The solution of \eqref{two-comp} continuously depends on initial 
	data, i.e., for any sequences $u_{0,n}(x)$, $v_{0,n}(x)$ such 
	that $u_{0,n}\to u_0$, $v_{0,n}\to v_0$ in $H^s$, $s>\frac{5}{2}$, as $n\to\infty$ we have
	\begin{equation}\label{cont-dep-lim}
	\lim\limits_{n\to\infty}
	\left(\|u-u_n\|
	_{C(I_{\delta_0}, H^{s})\cap C^1(I_{\delta_0}, H^{s-1})}
	+\|v-v_n\|
	_{C(I_{\delta_0}, H^{s})\cap C^1(I_{\delta_0}, H^{s-1})}
	\right)=0,
	\end{equation}
	with $I_{\delta_0}=[-T_{\delta_0},T_{\delta_0}]$.
	Here the pairs $u(x,t), v(x,t)$ and $u_n(x,t)$, $v_n(x,t)$ are the 
	solutions of the Cauchy problem for \eqref{two-comp} with initial data $u(x,0)=u_0(x)$, 
	$v(x,0)=v_0(x)$ and
	$u_n(x,0)=u_{0,n}(x)$, $v_n(x,0)=v_{0,n}(x)$, respectively.
\end{theorem}

\begin{remark}[Nonuniform continuity of the data-to-solution map]
	Observe that the data-to-solution map of the Cauchy problem 
	for \eqref{two-comp} is continuous in $C(I_{\delta_0}, H^s)$, $s>\frac{5}{2}$, 
	but not uniformly continuous.
	To show this consider \eqref{two-comp} with $v=u$.
	In such a case the two-component system reduces to the FORQ equation, 
	where the solution map is not uniformly continuous \cite{HM14}.
	Moreover, since the sequences of approximate solutions defined 
	in \cite[equations (5.1), (6.2)]{HM14}, satisfy the \textit{PT}-symmetry 
	condition (in (6.2) one should choose cutoff functions), the 
	data-to-solution map of the Cauchy problem for the nonlocal 
	FORQ equation \eqref{nFORQ-eq} is not uniformly continuous as well.
\end{remark}

Our work for the two-component system \eqref{two-comp} 
presents an alternative approach for solving the 
Cauchy problem. Theorems \ref{ex-un-sol} and \ref{cont-dep} 
correspond to specific cases of a broader result presented 
in \cite{MM13} (also see \cite{YQZ15}), with 
our method differing by utilizing a technique that reduces 
the problem to an ordinary differential equation (ODE) 
in a Banach space, as outlined in \cite{HM14}. 
We believe that our approach may prove useful for establishing 
global solutions, even in cases where singularities occur within finite time, 
as discussed in \cite{BC07, HR07} for the Camassa-Holm equation. 
We will come back to this in a future work.

In addition, in the next theorem, we establish new H{\"o}lder continuity 
properties of the data-to-solution map for \eqref{two-comp}. 
To the best of our knowledge, these properties 
have not been previously reported. Our approach builds upon 
the work presented in \cite{HM14-1} and involves refining 
estimates for specific values of the parameters $s$ and $r$, as 
indicated in \eqref{gamma} below. By employing these refined estimates, we have 
achieved improvements over the results obtained in \cite{HM14-1} for the FORQ equation. 
This improvement is demonstrated in Corollary \ref{FORQ-Holder}, where 
we showcase the enhanced H{\"o}lder continuity properties.

\begin{theorem}
	[H\"older continuity of data-to-solution map]
	\label{Holder-cont}
	Fix any $\rho>0$ and consider two solutions 
	$(u_j,v_j)$, $j=1,2$, of the 
	Cauchy problems for system \eqref{two-comp} with initial data
	$u_{0,j}=u_j(x,0)$ and
	$v_{0,j}=v_j(x,0)$, such that 
	$\|u_{0,j}\|_{H^{s}},\|v_{0,j}\|_{H^s}
	\leq\rho$, $s>\frac{5}{2}$, $j=1,2$.
	Then the two solutions
	$u_j,v_j\in
	C([-T_{\rho,\delta_0},T_{\rho,\delta_0}], H^{s})
	\cap 
	C^1([-T_{\rho,\delta_0},T_{\rho,\delta_0}], H^{s-1})$
	with 
	$T_{\rho,\delta_0}=
	\frac{1-\delta_0}{32C_s\rho^2}$, $0<\delta_0<1$,
	satisfy the following H\"older stability estimates:
	\begin{subequations}
		\begin{align}
		\label{H-cont-sol}
		&\|(u_1,v_1)(t)-(u_2,v_2)(t)\|_{H^{r}}
		\lesssim_{r,s,\rho,\delta_0}
		\|(u_{0,1},v_{0,1})
		-(u_{0,2},v_{0,2})\|_{H^{r}}^\gamma,&&
		t\in[-T_{\rho,\delta_0},T_{\rho,\delta_0}],\\
		\label{H-cont-pd-t}
		&\|\partial_t(u_1,v_1)(t)
		-\partial_t(u_2,v_2)(t)\|_{H^{p}}
		\lesssim_{r,s,\rho,\delta_0}
		\|(u_{0,1},v_{0,1})
		-(u_{0,2},v_{0,2})\|_{H^{p+1}}^\mu,&&
		t\in[-T_{\rho,\delta_0},T_{\rho,\delta_0}],
		\end{align}
	\end{subequations}
	where $\gamma\in(0,1]$ and 
	$\mu\in(0,1]$ has the following form, depending on 
	$(s,r)$ and $(s,p)$ respectively:
	\begin{equation}\label{gamma}
	\gamma=
	\begin{cases}
	1, &(s,r)\in A_1,\\
	\frac{2s-3}{s-r}, &(s,r)\in A_2,\\
	\frac{2(s-r)}{3+\ve_0}, &(s,r)\in A_3,\\
	\frac{2s-3-\ve_0}{2(s-r)}, &(s,r)\in A_4,\\
	s-r, &(s,r)\in A_5,\\
	\frac{s}{s-r}, &(s,r)\in A_6
	\end{cases}
	\qquad
	\mu=
	\begin{cases}
	1, &(s,p)\in B_1,\\
	\frac{2(s-p-1)}{3+\ve_1}, &(s,p)\in B_2,\\
	\frac{2s-3-\ve_1}{2(s-p-1)}, &(s,p)\in B_3,\\
	\frac{2s-4}{s-p-1}, &(s,p)\in B_4,\\
	s-p-1, &(s,p)\in B_5,\\
	\frac{s-1}{s-p-1}, &(s,p)\in B_6,
	\end{cases}
	\end{equation}
	with small $\ve_0,\ve_1>0$ such that
	$\ve_0\in(0,2s-5)$ for $s\in(5/2,3]$,
	$\ve_1\in(0,2s-5)$ for $s\in(5/2,11/4)$ and
	$\ve_1\in(0,2s-11/2)$ for $s\in(11/4,3)$.
	Here
	(see Figures \ref{A-regions} and \ref{B-regions}; cf.\,\,\cite[Figure 1]{HM14-1})
	\begin{equation}\label{A-j-regions}
	\begin{split}
	&A_1=\{(s,r):s>5/2,\,
	r\in(3/2,s-1]\}\cup
	\{(s,r):s>5/2,\,
	r\in[0,3/2],\\
	&\qquad\quad r\in[3-s,s-3/2)\},\\
	&A_2=\{(s,r):s\in(5/2,3),\,r<3-s\},\\
	&A_3=\{(s,r):s\in(5/2,3],\,
	r\in[s-3/2,3/2],\,
	4(s-r)^2>3(2s-3)\},\\
	&A_4=\{(s,r):s\in(5/2,3],\,r\in[s-3/2,3/2],\,
	4(s-r)^2\leq 3(2s-3)\},\\
	&A_5=\{(s,r):s>5/2,\,r\in(s-1,s)\},\\
	&A_6=\{(s,r): s\geq 3,\,r<0\},
	\end{split}
	\end{equation}
	and
	\begin{equation}
	\nonumber
	\begin{split}
	&B_1=\{(s,p):s>5/2,\,
	p\in(1/2,s-2]\}\cup
	\{(s,p):s>11/4,\,
	p\in[0,1/2],\,p\in[3-s,s-5/2)\},\\
	&B_2=\{(s,p):s\in(11/4,3),\,p\in[s-5/2,1/2],\,
	4(s-p-1)^2> 3(2s-3)\},\\
	&B_3=\{(s,p):s\in(5/2,11/4],\,p\leq 1/2\}
	\cup\{(s,p):s\in(11/4,3],\,
	p\in[s-5/2,1/2],\\
	&\qquad\quad 4(s-p-1)^2\leq 3(2s-3)\},\\
	&B_4=\{(s,p):s\in(11/4,3),\,p<3-s\},\\
	&B_5=\{(s,p):s>5/2,\,p\in(s-2,s-1)\},\\
	&B_6=\{(s,p):s\geq 3,\,p<0\}.
	\end{split}
	\end{equation}
	\begin{figure}
		\centering{\includegraphics[scale=0.5]
			{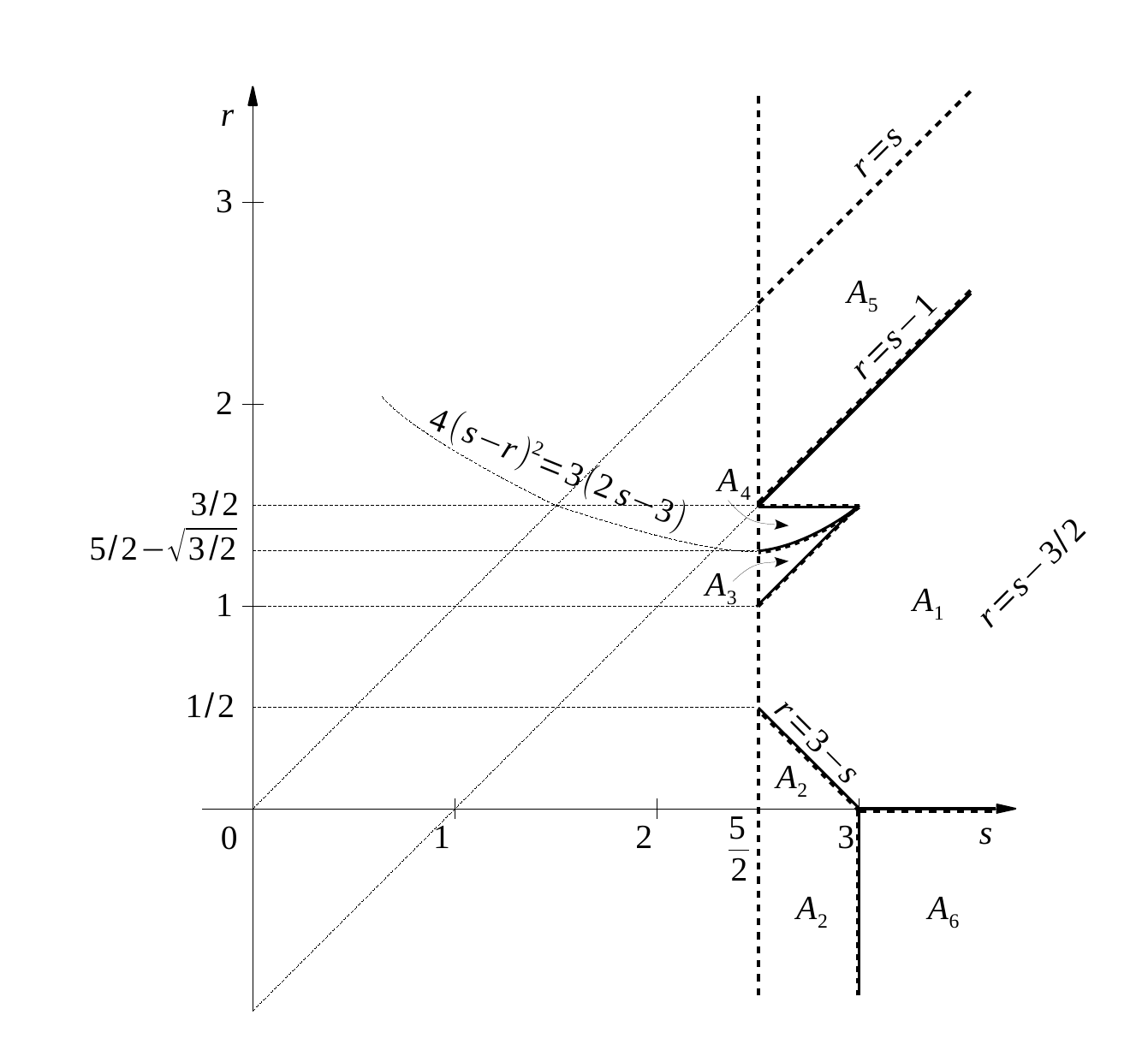}}
		\caption{Regions $A_j$, $j=1,\dots,6$ in the $(s,r)$ plane.}
		\label{A-regions}
	\end{figure}
	\begin{figure}
		\centering{\includegraphics[scale=0.5]
			{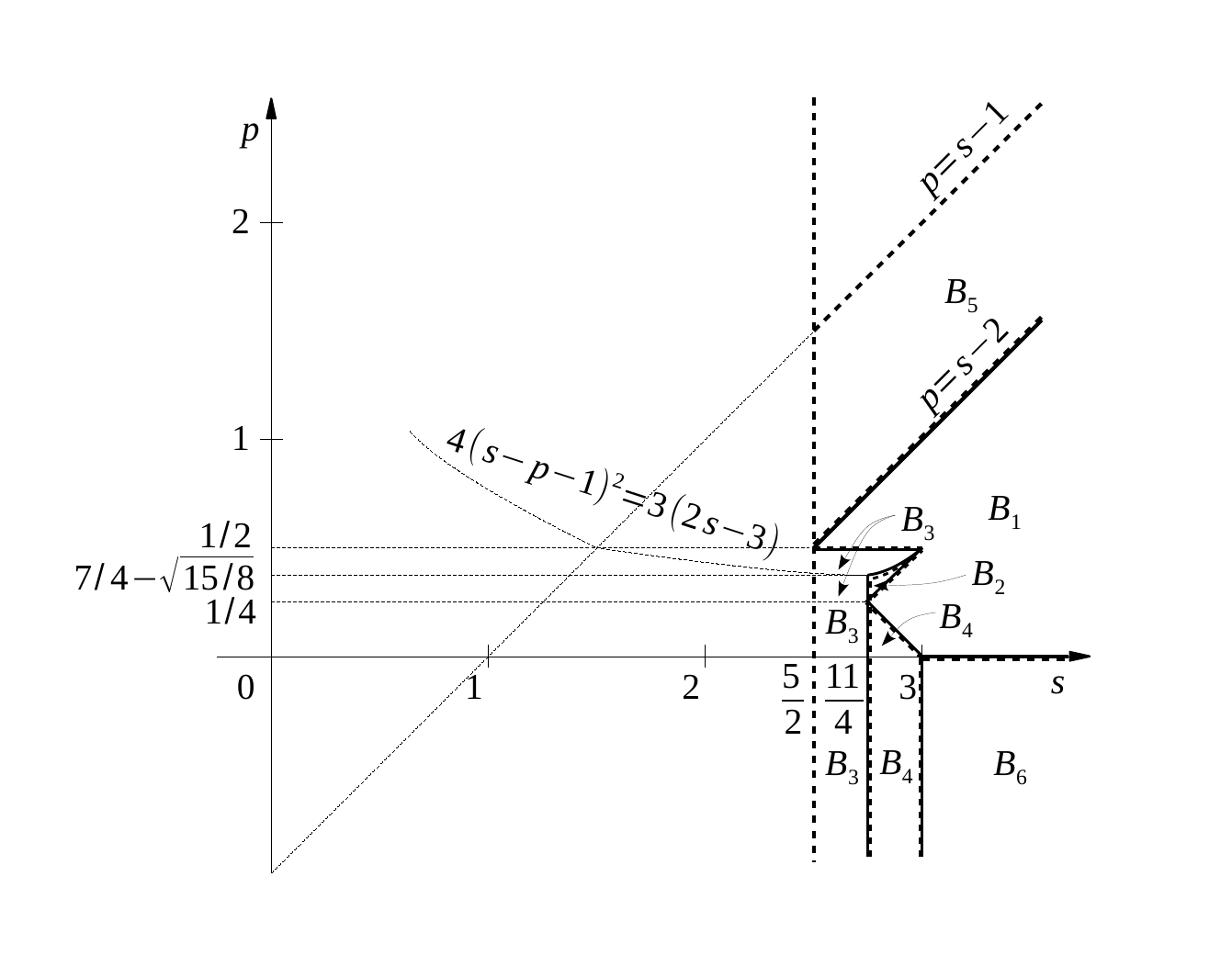}}
		\caption{Regions $B_j$, $j=1,\dots,6$ in the $(s,p)$ plane.}
		\label{B-regions}
	\end{figure}
\end{theorem}

\begin{remark}
	The function $\gamma=\gamma(s,r)$ is continuous
	along the lines $\{(s,r):s=3,\,r<0\}$,
	$\{(s,r):s\geq3,\,r=0\}$,
	$\{(s,r):s\in(5/2,3],\,r=3-s\}$
	and
	$\{(s,r):s>5/2,\,r=s-1\}$
	as well as,
	in the limit $\ve_0\to0$, along the lines 
	$\{(s,r):s\in(5/2,3],\,r=s-3/2\}$,
	$\{(s,r):s\in(5/2,3],\,r=3/2\}$ and along the parabola
	$4(s-r)^2=3(2s-3)$ for $s\in(5/2,3]$.
	
	Concerning the function $\mu=\mu(s,p)$, it is 
	continuous
	along the lines $\{(s,p):s=3,\,p<0\}$,
	$\{(s,p):s\geq3,\,p=0\}$,
	$\{(s,p):s\in(11/4,3],\,p=3-s\}$,
	$\{(s,p):s>5/2,\,p=s-2\}$
	as well as,
	in the limit $\ve_1\to0$, along the lines
	$\{(s,p):s\in(11/4,3],\,p=s-5/2\}$,
	$\{(s,p):s\in(5/2,3],\,p=1/2\}$ and the parabola
	$4(s-p-1)^2=3(2s-3)$ for $s\in(11/4,3]$.
\end{remark}

The main focus of our work revolves around 
Theorems \ref{ex-un-sol}, \ref{cont-dep}, and \ref{Holder-cont}, where 
we present our key findings regarding the local well-posedness of the 
Cauchy problem associated with the two-component peakon 
system \eqref{two-comp}. In order to provide a clearer perspective, let us 
emphasize some specific instances of the aforementioned theorems 
that correspond to reductions of \eqref{two-comp}:

\begin{corollary}
Taking $v(x,t)=u(-x,-t)$ in Theorems \ref{ex-un-sol}, \ref{cont-dep} and \ref{Holder-cont}, 
we obtain local existence, uniqueness and continuous dependence 
on the initial data in $C([-T_{\delta_0},T_{\delta_0}], H^{s})
\cap C^1([-T_{\delta_0},T_{\delta_0}], H^{s-1})$, $s>\frac{5}{2}$, 
as well as the H\"older continuity properties of the solution of 
the Cauchy problem for the nonlocal FORQ equation \eqref{nFORQ-eq}.
\end{corollary}

\begin{corollary}\label{FORQ-Holder}
Taking $v(x,t)=u(x,t)$ in Theorem \ref{Holder-cont}, we obtain the 
improved H\"older exponent in the regions $A_1$, $A_3$ and $A_4$ 
comparing with that in \cite{HM14-1} for the FORQ 
equation (see Remarks \ref{Lip-A-1} and \ref{Hold-exp-3-4} below for details). 
Also, notice that here we present the H\"older continuity 
properties in $H^r$ for all $r<s$, $s>\frac{5}{2}$.
\end{corollary}

The article is structured as follows. We begin in 
Section \ref{Prelim} by introducing relevant notations and gathering 
key mathematical results that will be utilized throughout the paper. 
In Section \ref{Hamilt-dual}, we explore the relationship 
between the two-component system \eqref{two-comp} and 
a specific member of the AKNS hierarchy, thereby providing some 
background information. Next, in Section \ref{loc-ex-un}, we establish the 
local existence and uniqueness of the solution, leading us to 
Theorem \ref{ex-un-sol}. Subsequently, in Section \ref{Cont}, we 
focus on proving Theorem \ref{cont-dep} by demonstrating 
the continuous dependence of the solution on the initial data. 
Finally, we delve into the investigation of 
the H{\"o}lder continuity properties of 
Theorem \ref{Holder-cont} in Section \ref{Holder}.

\section{Preliminaries}\label{Prelim}
In this section we introduce some notations as well as 
important mathematical results to be used throughout the paper.

The Fourier transform of the function $f(x)$, $x\in A$ has the form
$$
\mathcal{F}(f)(k)=
\frac{1}{2\pi}\int_Ae^{-\mathrm{i}kx}f(x)\,dx,\quad k\in\mathbb{R},
$$
with the inverse relation given by
$$
\mathcal{F}^{-1}(F)(x)=
\int_\mathbb{R}e^{\mathrm{i}kx}F(k)\,dk,\quad\mbox{for\,\,} x\in\mathbb{R},\quad
\mathcal{F}^{-1}(F)(x)=
\sum\limits_{k\in\mathbb{Z}}e^{\mathrm{i}kx}F(k)
,\quad\mbox{for\,\,} x\in\mathbb{R}/2\pi\mathbb{Z}.
$$
In these notations the Plancherel formula is as follows
\begin{equation}\label{Planch}
(f,g)_{L^2(A)}=2\pi(\mathcal{F}(f),\mathcal{F}(g))_{L^2(A)},
\end{equation}
and
\begin{equation}\label{conv}
\mathcal{F}(f\ast g)=2\pi\mathcal{F}(f)\mathcal{F}(g).
\end{equation}
Introduce the smooth function
$$
j(x)=
\begin{cases}
I_1\exp\left\{\frac{1}{x^2-1}\right\},&|x|<1,\\
0& |x|\geq1,
\end{cases}
$$
with $I_1=\left(\int_{-1}^{1}
\exp\left\{\frac{1}{x^2-1}\right\}\,dx\right)^{-1}$.
Then define a Friedrichs mollifier $J_\ve$ of the function $f$ as follows:
$$
J_\ve f\equiv f^\ve=j_\ve\ast f,
$$
where
$j_\ve(x)=\frac{1}{\ve}j\left(\frac{x}{\ve}\right)$ in the case $A=\mathbb{R}$
and 
$j_\ve(x)=\sum\limits_{k\in\mathbb{Z}}
\mathcal{F}(j)(\ve k)e^{\mathrm{i}kx}$ in the periodic case.
Notice that in view of \eqref{Planch} and \eqref{conv} we have
\begin{equation}\label{moll-inner-prod}
(J_\ve f,g)_{L^2(A)}=(f,J_\ve g)_{L^2(A)}.
\end{equation}
The norm of the function $f(x)$ in Sobolev 
space $H^s(A)$, $s\in\mathbb{R}$ is defined by
$$
\|f\|_{H^s(A)}^2=\|D^sf\|_{L^2(A)}^2=
2\pi\int_A(1+|k|^2)^s|\mathcal{F}(f)(k)|^2\,dk,
$$
where 
$$
D^s=(1-\partial_x^2)^\frac{s}{2},\quad s\in\mathbb{R},
$$ 
is the Bessel potential of order $s$.

To enhance readability, we will adopt the notation $H^s$ instead 
of explicitly writing $H^s(A)$ throughout the subsequent sections, 
as long as it does not cause confusion for the reader.

We will leverage a set of commutator estimates, which are 
succinctly summarized in the following lemma:
\begin{lemma}
	Let $[f,g]=fg-gf$ be a commutator.
	Then the following estimates hold:
	\begin{enumerate}
	\item the Kato-Ponce inequality \cite{KP88}:
	\begin{equation}\label{Kato-Ponce}
	\|[D^s,f]g\|_{L^2(A)}\lesssim_s
	(\|D^sf\|_{L^2(A)}\|g\|_{L^\infty(A)}
	+\|\partial_xf\|_{L^\infty(A)}\|D^{s-1}g\|_{L^2(A)}),\quad s>0,
	\end{equation}
	with $f,g\in H^s(A)$ and 
	$\partial_xf,g\in L^\infty(A)$;

	\item the Calderon-Coifman-Meyer commutator 
	estimate (see, e.g., \cite[Proposition 4.2]{T03}):
	\begin{equation}
	\label{CCM}
	\|[D^\sigma\partial_x,f]g\|_{L^2(A)}
	\lesssim_{\nu,\sigma}
	\|f\|_{H^\nu(A)}
	\|g\|_{H^\sigma(A)},
	\end{equation}
	provided $f\in H^{\nu}$, $g\in H^\sigma$, $\nu>\frac{3}{2}$ and $0\leq\sigma+1\leq\nu$;

	\item an inequality of Himonas and Kenig (see \cite[Lemma 2]{HK09}):
	\begin{equation}\label{HK}
	\|[J_\ve,f]\partial_xg\|_{L^2(A)}\lesssim
	\|\partial_x f\|_{L^\infty(A)}\|g\|_{L^2(A)},
	\end{equation}
	where $\partial_xf\in L^\infty(A)$ and $g\in L^2(A)$;
	
	\item an inequality of Himonas and Holmes (see \cite[Lemma 3]{HH13}):
	\begin{equation}
	\label{HH}
	\|fg\|_{H^{\hat{r}-1}}
	\lesssim_{\hat{r},\hat{s}}
	\|f\|_{H^{\hat{s}-1}}
	\|g\|_{H^{\hat{r}-1}},
	\end{equation}
	where $\hat{r}\in[0,1]$, $\hat{s}>3/2$ and $\hat{r}+\hat{s}\geq 2$.
	\end{enumerate}
\end{lemma}
We will use the Sobolev embedding theorem
\begin{equation}\label{Sob-emb}
\|f\|_{L^{\infty}(A)}\lesssim_s
\|f\|_{H^s(A)},\quad s>\frac{1}{2},
\end{equation}
the algebra property of the Sobolev space $H^s$
\begin{equation}\label{alg}
\|fg\|_{H^s(A)}\lesssim_s
\|f\|_{H^s(A)}\|g\|_{H^s(A)},\quad s>\frac{1}{2},
\end{equation}
as well as
\begin{equation}\label{ineq}
\|D^mf\|_{H^s(A)}\leq\|f\|_{H^{s+m}(A)},\,\, \|\partial_xf\|_{H^s(A)}\leq\|f\|_{H^{s+1}(A)},\,\,
\|w_\ve^\ve\|_{H^s(A)}\leq\|w_\ve\|_{H^s(A)}.
\end{equation}
Also we will employ the following interpolation inequality
(see, e.g., \cite[Lemma 5]{HM14-1})
\begin{equation}\label{interpolation}
\|f\|_{H^\sigma}\leq
\|f\|_{H^{\sigma_1}}
^{\frac{\sigma_2-\sigma}{\sigma_2-\sigma_1}}
\|f\|_{H^{\sigma_2}}
^{\frac{\sigma-\sigma_1}{\sigma_2-\sigma_1}},
\end{equation}
provided $\sigma_1<\sigma<\sigma_2$ and $f\in H^\sigma$.

Finally, will apply the following variant of the Rellich-Kondrachov theorem 
\begin{lemma}[Rellich-Kondrachov]
	\label{RK}
	Let $W\subset\mathbb{R}$ be a compact set.
	Then the space $H^{s'}(W)$ is compactly 
	embedded in $H^{s''}(W)$ for any $0<s''<s'<\infty$:
	$$
	H^{s'}(W)\Subset H^{s''}(W).
	$$
\end{lemma}

\begin{proof}
The proof follows by utilizing the following three key components: 
a) the fact that $H^s(W)\Subset L^2(W)$ holds for all $s\geq 1$, 
b) citing \cite[Theorem 2, ch. 1.16.4]{T78}, and 
c) referring to \cite[ch. 2.4.1]{T78}.
\end{proof}


\section{Peakon system (\ref{two-comp}) and AKNS hierarchy}
\label{Hamilt-dual}

In \cite{OR96} the authors proposed a tool for the systematic 
derivation of new integrable equations with nonlinear dispersion.
According to their method, one considers a three parameter family of operators, say
$J=a_1K_1+a_2K_2+a_3K_3$, which is Hamiltonian (i.e., it is skew-symmetric 
and satisfies the Jacobi identity) for any $a_j\in\mathbb{R}$, $j=1,2,3$.
Then define four operators $J_j$, $\tilde{J}_j$, $j=1,2$ as particular 
reductions of $J$ (for example, one can take $(J_1, J_2)=(K_1,K_2+K_3)$
and $(\tilde J_1,\tilde J_2)=(K_1+K_2,K_3)$).
Notice that since $J$ is Hamiltonian, the pairs $(J_1, J_2)$ 
and $(\tilde J_1, \tilde J_2)$ are compatible, meaning that 
both operators $(\alpha_1J_1+\alpha_2J_2)$ and
$(\tilde\alpha_1\tilde J_1+\tilde\alpha_2\tilde J_2)$ 
are Hamiltonian for any 
$\alpha_j,\tilde\alpha_j\in\mathbb{R}$,
$j=1,2$.
Thus the pair of Hamiltonian operators
$(J_1, J_2)$
generates the first bi-Hamiltonian hierarchy, while the other pair
$(\tilde J_1,\tilde J_2)$ generates the dual bi-Hamiltonian hierarchy.

In \cite[Example 2]{OR96} (see also \cite{KLOQ16}) it was 
shown that the FORQ equation is dual to the mKdV equation
\begin{equation}\label{mkdv}
\partial_t\tilde m
+\partial_x^3\tilde m
+6\tilde m^2\partial_x\tilde m=0,\quad
\tilde m=\tilde m(x,t).
\end{equation}
Therefore, a natural inclination arises to establish a 
duality relationship between the two-component equation \eqref{two-comp} 
and its counterpart in the AKNS hierarchy.
Notably, \cite{TL13} successfully established the duality 
between \eqref{two-comp} and the member 
belonging to the Wadati-Konno-Ichikawa hierarchy \cite{WKI79}.

Here we consider, in contrast to \cite{OR96}, \textit{two} three-parameter 
families of Hamiltonian operators $L_1$ and $L_2$ (see \eqref{L_1-L_2} below) 
and obtain from them two pairs of compatible Hamiltonian operators $(J_1,J_2)$ and
$(\tilde J_1, \tilde J_2)$ (see \eqref{J_1} and \eqref{J_2} below) 
by choosing different values of the parameters. 
The pair $(J_1, J_2)$ generates the hierarchy of the 
system \eqref{two-comp} (cf.~\cite{SQQ11, TL13, XQZ15}), while 
the pair $(\tilde J_1, \tilde J_2)$ gives rise to the 
AKNS hierarchy (cf.~\cite{T89} and, e.g., \cite[Subsection 2.1]{LM21}).

Consider two operators
\begin{equation}
	\label{L_1-L_2}
	\begin{split}
	&L_1=a_1
	\begin{pmatrix}
	0&1\\
	-1&0
	\end{pmatrix}
	+b_1
	\begin{pmatrix}
	0&\partial_x\\
	\partial_x&0
	\end{pmatrix}
	+c_1
	\begin{pmatrix}
	0&\partial_x^2\\
	-\partial_x^2&0
	\end{pmatrix},\\
	&L_2=a_2
	\begin{pmatrix}
	0&\partial_x\\
	\partial_x&0
	\end{pmatrix}
	+b_2
	\begin{pmatrix}
	-m\partial_x^{-1}m&m\partial_x^{-1}n\\
	n\partial_x^{-1}m&-n\partial_x^{-1}n
	\end{pmatrix}
	+c_2
	\begin{pmatrix}
	\partial_xm\partial_x^{-1}m\partial_x
	&\partial_xm\partial_x^{-1}n\partial_x\\
	\partial_xn\partial_x^{-1}m\partial_x
	&\partial_xn\partial_x^{-1}n\partial_x
	\end{pmatrix},
	\end{split}
\end{equation}
with $a_j,b_j,c_j\in\mathbb{R}$, $j=1,2$.
From \cite[Proposition 1]{XQZ15} it follows that both $L_1$ and $L_2$ are Hamiltonian operators.
Taking  
$(a_1,b_1,c_1)=(0,1,1)$ 
and
$(a_1,b_1,c_1)=(1,0,0)$ 
we obtain operators $J_1$ and $\tilde{J}_1$ respectively, which have the form
\begin{equation}\label{J_1}
J_1=
\begin{pmatrix}
0&\partial_x+\partial_x^2\\
\partial_x-\partial_x^2&0
\end{pmatrix},\quad
\tilde J_1=
\begin{pmatrix}
0&1\\
-1&0
\end{pmatrix},
\end{equation}
while for 
$(a_2,b_2,c_2)=(0,0,1)$ and
$(a_2,b_2,c_2)=(-1,-2,0)$, we have operators $J_2$ and $\tilde{J}_2$ respectively (for the latter we use 
$(\tilde m,\tilde n)$ instead of 
$(m,n)$):
\begin{equation}\label{J_2}
J_2=
\begin{pmatrix}
\partial_xm\partial_x^{-1}m\partial_x
&\partial_xm\partial_x^{-1}n\partial_x\\
\partial_xn\partial_x^{-1}m\partial_x
&\partial_xn\partial_x^{-1}n\partial_x
\end{pmatrix},\quad
\tilde J_2=
\begin{pmatrix}
2\tilde m\partial_x^{-1}\tilde m
&-\partial_x
-2\tilde m\partial_x^{-1}\tilde n\\
-\partial_x-2\tilde n\partial_x^{-1}\tilde m
&2\tilde n\partial_x^{-1}\tilde n
\end{pmatrix}.
\end{equation}
In view of \cite[Proposition 1]{XQZ15}, the pairs of operators $(J_1,J_2)$ and $(\tilde J_1,\tilde J_2)$ are compatible.
Therefore Magri's theorem \cite{M78} (see also \cite[Proposition 2]{FF81}) implies that there exist two infinite hierarchies of bi-Hamiltonian systems:
\begin{equation}\label{bi-H-two-comp}
\partial_t(m,n)^T
=(J_2J_1^{-1})^i[\partial_x(m,n)^T]
=J_1\left(\frac{\delta H_{i+1}}{\delta m},
\frac{\delta H_{i+1}}{\delta n}\right)^T
=J_2\left(\frac{\delta H_{i}}{\delta m},
\frac{\delta H_{i}}{\delta n}\right)^T, 
\end{equation}
for $i=0,1,2,\dots$ and
\begin{equation}\label{bi-H-mkdv}
\partial_t(\tilde m,\tilde n)^T=
(\tilde J_2\tilde J_1^{-1})^j
[\partial_x(\tilde m,\tilde n)^T]
=\tilde J_1\left(\frac
{\delta\tilde H_{j+1}}{\delta\tilde m},
\frac
{\delta\tilde H_{j+1}}
{\delta\tilde n}\right)^T
=\tilde J_2\left(\frac
{\delta\tilde H_{j}}{\delta\tilde m},
\frac
{\delta\tilde H_{j}}
{\delta\tilde n}\right)^T,
\end{equation}
for $j=0,1,2,\dots$.

The bi-Hamiltonian systems \eqref{bi-H-two-comp} form 
the hierarchy of \eqref{two-comp}, where the latter can be 
obtained for $i=1$ (cf.~\cite[Section IV]{TL13} 
and \cite[(55)--(57)]{XQZ15}):
\begin{equation}
\nonumber
\eqref{two-comp}\Leftrightarrow
\partial_t(m,n)^T=
J_1\left(\frac{\delta H_{2}}{\delta m},
\frac{\delta H_{2}}{\delta n}\right)^T
=J_2\left(\frac{\delta H_{1}}{\delta m},
\frac{\delta H_{1}}{\delta n}\right)^T,
\end{equation}
with 
\begin{equation}
\nonumber
H_1=\int_Am(v+\partial_xv)\,dx,\quad
H_2=\frac{1}{2}\int_A(u-\partial_xu)^2
(v+\partial_xv)n\,dx.
\end{equation}
Equations \eqref{bi-H-mkdv}, on the other hand, belong 
to the set of equations known as the 
AKNS hierarchy, as discussed in \cite[Subsection 2.1]{LM21}.

Notice that the system 
\begin{equation}\label{mkdvh}
	\begin{split}
	&\partial_t \tilde m
	+\partial_x^3\tilde m
	+6\tilde m\tilde n\partial_x\tilde m=0,
	\quad \tilde m=\tilde m(x,t)
	,\,\tilde n=\tilde n(x,t),\\
	&\partial_t \tilde n
	+\partial_x^3\tilde n
	+6\tilde n\tilde m\partial_x\tilde n=0,
	\end{split}
\end{equation}
which reduces to the mKdV equation \eqref{mkdv} in the 
case $\tilde m=\tilde n$ and to the nonlocal mKdV \eqref{nmkdv} 
in the case $\tilde n(x,t)=\tilde m(-x,-t)$, can be 
obtained from \eqref{bi-H-mkdv} with $j=2$ \cite{T89}:
\begin{equation}
\nonumber
\eqref{mkdvh}\Leftrightarrow
\partial_t(\tilde m,\tilde n)^T=
\tilde J_1\left(\frac
{\delta\tilde H_{3}}{\delta\tilde m},
\frac
{\delta\tilde H_{3}}
{\delta\tilde n}\right)^T
=\tilde J_2\left(\frac
{\delta\tilde H_{2}}{\delta\tilde m},
\frac
{\delta\tilde H_{2}}
{\delta\tilde n}\right)^T,
\end{equation}
where 
\begin{equation}
\nonumber
\tilde H_{2}=\int_A 
(\tilde m^2\tilde n^2
-(\partial_x\tilde m)
\partial_x\tilde n)\,dx,\quad
\tilde H_{3}=\int_A 
(3\tilde m^2\tilde n\partial_x\tilde n
-(\partial_x\tilde m)
\partial_x^2\tilde n)\,dx.
\end{equation}

Thus we arrive at the following
\begin{proposition}
Consider the three-parameter families of Hamiltonian operators 
$L_1(a_1,b_1,c_1)$ and 
$L_2(a_2,b_2,c_2)$ defined in \eqref{L_1-L_2}.
Then the Hamiltonian operators 
(see \eqref{J_1} and \eqref{J_2})
$$
J_1=L_1(0,1,1),\,\,
\tilde{J}_1=L_1(1,0,0),\,\,
J_2=L_2(0,0,1),\,\
\tilde{J}_2=L_2(-1,-2,0)
$$
satisfy the following properties:
\begin{enumerate}
\item the pairs $(J_1, J_2)$ and $(\tilde{J}_1, \tilde{J}_2)$ are 
compatible pairs of operators,
\item the recursion operators $J_2J_1^{-1}$ and $\tilde{J}_2\tilde{J}_1^{-1}$ 
generate the infinite hierarchies of bi-Hamiltonian systems 
\eqref{bi-H-two-comp} and \eqref{bi-H-mkdv}, respectively,
where the former reduces to \eqref{two-comp} for $i=1$, 
while the latter gives rise to the member of 
AKNS hierarchy \eqref{mkdvh} for $j=2$.
\end{enumerate}
\end{proposition}

\section{Local existence and uniqueness}\label{loc-ex-un}

To investigate the local well-posedness issues related 
to the Cauchy problem concerning the two-component equation \eqref{two-comp}, 
we employ a reduction technique that transforms it into 
an ODE system in a Banach space. This reduction allows us 
to address the problem more effectively, as demonstrated 
in \eqref{moll-ODE} (see also \cite{HM14}). Before considering the 
appropriate ODE counterpart, we establish an equivalence 
between \eqref{two-comp} and a system of four PDEs. 
This equivalence will facilitate our analysis.

\begin{proposition}\label{two-comp-to-syst}
	The pair $(u,v)$ solves the Cauchy problem 
	for the two-component equation \eqref{two-comp} in 
	$C\left([-T,T], \left(H^{s}(A)\right)^2\right)$, $s>\frac{5}{2}$, 
	with initial data $u(x,0)=u_0(x)$ and $v(x,0)=v_0(x)$ if and only if 
	the function
	$(u,w,v,z)\equiv(u,\partial_xu,v,\partial_xv)$ 
	satisfies the following system in $C\left([-T,T], \left(H^{s-1}(A)\right)^4\right)$:
	\begin{subequations}
	\label{ODE}
	\begin{align}
	\label{ODE-a}
	&\partial_t u =-\frac{1}{3}w^2z+\frac{1}{3}(2uwv+u^2z)
	+F(u,w,v,z),\\
	\label{ODE-b}
	\nonumber
	&\partial_t w=-w(\partial_xw)z-\frac{1}{3}w^2v+\frac{4}{3}uwz
	+\frac{1}{3}(2u(\partial_xw)v+u^2\partial_xz)
	-\frac{2}{3}u^2v\\
	&\qquad\quad+\{u(\partial_xw)z-w(\partial_xw)v+uwv-u^2z\}
	+G(u,w,v,z),\\
	\label{ODE-c}
	&\partial_tv=-\frac{1}{3}wz^2+\frac{1}{3}(2uvz+wv^2)
	+\hat F(u,w,v,z),\\
	\label{ODE-d}
	\nonumber
	&\partial_tz
	=-wz\partial_xz-\frac{1}{3}uz^2+\frac{4}{3}wvz
	+\frac{1}{3}(2uv\partial_xz+(\partial_xw)v^2)
	-\frac{2}{3}uv^2\\
	&\qquad\quad-\{wv\partial_xz-uz\partial_xz+uvz-wv^2\}
	+\hat G(u,w,v,z),\\
	&
	u(x,0)=u_0(x),\quad w(x,0)=w_0(x),\quad
	v(x,0)=v_0(x),\quad
	z(x,0)=z_0(x),
	\end{align}
	\end{subequations}
	where 
	$w_0=\partial_xu_0$, $z_0=\partial_xv_0$ and
	\begin{equation}
	\label{nonl-terms}
	\begin{split}
	&F=D^{-2}\left(
	\frac{1}{3}w^2z+
	\left\{uw\partial_xz-w(\partial_xw)v
	+\frac{1}{3}u(u\partial_x^2z-(\partial_x^2w)v)
	\right\}
	\right)
	+D^{-2}\partial_x\left(
	\frac{2}{3}u^2v+w^2v+B
	\right),\\
	&\hat F=D^{-2}\left(
	\frac{1}{3}wz^2-
	\left\{uz\partial_xz-(\partial_xw)vz
	+\frac{1}{3}v(u\partial_x^2z-(\partial_x^2w)v)
	\right\}
	\right)
	+D^{-2}\partial_x\left(
	\frac{2}{3}uv^2+uz^2+\hat B
	\right),\\
	& G=D^{-2}\partial_x\left(
	\frac{1}{3}w^2z+
	\left\{uw\partial_xz-w(\partial_xw)v
	+\frac{1}{3}u(u\partial_x^2z-(\partial_x^2w)v)
	\right\}
	\right)
	+D^{-2}\left(
	\frac{2}{3}u^2v+w^2v+B
	\right),\\
	& \hat G=D^{-2}\partial_x\left(
	\frac{1}{3}wz^2-
	\left\{uz\partial_xz-(\partial_xw)vz
	+\frac{1}{3}v(u\partial_x^2z-(\partial_x^2w)v)
	\right\}
	\right)
	+D^{-2}\left(
	\frac{2}{3}uv^2+uz^2+\hat B
	\right),
	\end{split}
	\end{equation}
	with
	\begin{align*}
	&B(u,w,v,z)=-u(\partial_xw)z+w(\partial_xw)v-uwv+u^2z
	+\frac{1}{3}(w(\partial_xw)z-w^2\partial_xz),\\
	&\hat B(u,w,v,z)=
	wv\partial_xz-uz\partial_xz+uvz-wv^2
	+\frac{1}{3}(wz\partial_xz-(\partial_xw)z^2).
	\end{align*}
\end{proposition}
\begin{proof}
	Let $(u,v)$ solve the Cauchy problem for \eqref{two-comp}. 
	Taking into account that
	\begin{equation}
	\nonumber
	\begin{split}
	-\frac{1}{3}D^2((\partial_xu)^2\partial_xv)=&
	-\frac{1}{3}(\partial_xu)^2\partial_xv
	+\frac{2}{3}(\partial_x^2u)^2\partial_xv
	+\frac{4}{3}(\partial_xu)(\partial_x^2u)\partial_x^2v
	+\frac{2}{3}(\partial_xu)(\partial_x^3u)\partial_xv\\
	&+\frac{1}{3}(\partial_xu)^2\partial_x^3v,
	\end{split}
	\end{equation}
	and
	\begin{align}
	\nonumber
	&\frac{1}{3}D^2(2u(\partial_xu)v+u^2\partial_xv)=
	\frac{1}{3}(2u(\partial_xu)v+u^2\partial_xv)\\
	\nonumber
	&\qquad-2\left(\frac{1}{6}u^2\partial_x^3v+
	u(\partial_xu)\partial_x^2v+(\partial_xu)^2\partial_xv
	+u(\partial_x^2u)\partial_xv+(\partial_xu)(\partial_x^2u)v
	+\frac{1}{3}u(\partial_x^3u)v\right),
	\end{align}
	we conclude that \eqref{two-comp-a} can be 
	written in the form (compare with \cite[(2.1)]{HM14})
	\begin{equation}\label{nonl-nFORQ}
	\partial_tu=-\frac{1}{3}(\partial_xu)^2\partial_xv
	+\frac{1}{3}(2u(\partial_xu)v+u^2\partial_xv)+\Phi(u,v),
	\end{equation}
	with
	\begin{align}
	\nonumber
	\Phi(u,v)=&
	D^{-2}\left(
	\frac{1}{3}(\partial_xu)^2\partial_xv
	\right)
	+D^{-2}\partial_x\left(
	\frac{2}{3}u^2v+(\partial_xu)^2v+B
	\right)\\
	\nonumber
	&+D^{-2}\left(
	u(\partial_xu)\partial_x^2v-(\partial_xu)(\partial_x^2u)v
	+\frac{1}{3}u(u\partial_x^3v-(\partial_x^3u)v)
	\right).
	\end{align}
	Thus, we have verified \eqref{ODE-a}. 
	By taking the derivative of \eqref{nonl-nFORQ} with 
	respect to $x$ and combining similar terms, we can 
	deduce \eqref{ODE-b}.
	
	For \eqref{two-comp-b} we use
	\begin{equation}
	\nonumber
	\begin{split}
	-\frac{1}{3}D^2((\partial_xu)(\partial_xv)^2)=&\,
	\frac{2}{3}(\partial_xu)(\partial_xv)\partial_x^3v
	-\frac{1}{3}(\partial_xu)(\partial_xv)^2
	+\frac{2}{3}(\partial_xu)(\partial_x^2v)^2
	+\frac{4}{3}(\partial_x^2u)(\partial_xv)\partial_x^2v\\
	&+\frac{1}{3}(\partial_x^3u)(\partial_xv)^2,
	\end{split}
	\end{equation}
	and 
	\begin{align}
	\nonumber
	&\frac{1}{3}D^2(2uv\partial_xv+(\partial_xu)v^2)=
	\frac{1}{3}(2uv\partial_xv+\partial_xuv^2)\\
	\nonumber
	&\qquad-2\left(\frac{1}{3}uv\partial_x^3v+
	u(\partial_xv)\partial_x^2v+(\partial_xu)v\partial_x^2v
	+(\partial_xu)v\partial_x^2v+(\partial_x^2u)v(\partial_xv)
	+\frac{1}{6}(\partial_x^3u)v^2\right),
	\end{align}
	which implies that
	\begin{equation}\label{nonl-nFORQ-v}
	\partial_tv=-\frac{1}{3}(\partial_xu)(\partial_xv)^2
	+\frac{1}{3}(2uv(\partial_xv)+(\partial_xu)v^2)
	+\hat\Phi(u,v),
	\end{equation}
	with
	\begin{align}
	\nonumber
	\hat\Phi(u,v)=&
	D^{-2}\left(
	\frac{1}{3}(\partial_xu)(\partial_xv)^2
	\right)
	+D^{-2}\partial_x\left(
	\frac{2}{3}uv^2+u(\partial_xv)^2+\hat B
	\right)\\
	\nonumber
	&-D^{-2}\left(
	u(\partial_xv)\partial_x^2v-(\partial_x^2u)v(\partial_xv)
	+\frac{1}{3}v(u\partial_x^3v-(\partial_x^3u)v)
	\right),
	\end{align}
	which is equivalent to \eqref{ODE-c}.
	Differentiating \eqref{nonl-nFORQ-v} in $x$ we obtain \eqref{ODE-d}.
	
	For proving that the solution of \eqref{ODE} solves the Cauchy 
	problem for \eqref{two-comp} it is enough to show that from \eqref{ODE} it follows that
	that $w=\partial_xu$ and $z=\partial_xv$.
	Let $y_1=\partial_xu-w$ and $y_2=\partial_xv-z$, then direct 
	calculations show that $(y_1,y_2)$ satisfies the following 
	Cauchy problem for the linear ODE in the Banach 
	space $C\left([-T,T], \left(H^{s-1}(A)\right)^2\right)$:
	\begin{align*}
		&\frac{d}{dt}y_1=a_1y_1+a_2y_2,\\
		&\frac{d}{dt}y_2=a_3y_1+a_1y_2,\\
		&y_1(0)=0,\,\,y_2(0)=0,
	\end{align*}
	with $a_1=\frac{2}{3}(uz+wv)$, $a_2=\frac{2}{3}uw$ 
	and $a_3=\frac{2}{3}vz$, which has a unique solution $y_1=y_2=0$.
\end{proof}

\begin{remark}
	Observe that in the case $v=u$, $z=w$ the terms in curly brackets 
	and the functions $B$, $\hat B$ in \eqref{ODE}, \eqref{nonl-terms} vanish.
	Therefore the system \eqref{ODE} reduces to that 
	corresponding to the FORQ equation, cf.\,\,\cite[(2.4) with $-t$ 
	instead of $t$]{HM14}. Also notice that in the nonlocal 
	case (i.e., when $v(x,t)=u(-x,-t)$) the 
	equations \eqref{ODE-c}-\eqref{ODE-d} can be 
	obtained directly from \eqref{ODE-a}-\eqref{ODE-b}.
\end{remark}

\begin{remark}\label{u_x=w}
	Notice that in the proof of Proposition \ref{two-comp-to-syst} 
	we have shown that given a solution $(u,w,v,z)$ to \eqref{ODE} it 
	follows that $w=\partial_xu$ and $z=\partial_xv$.
\end{remark}
\begin{remark}
	A system similar to \eqref{ODE} was previously obtained in \cite[Section 3]{MM13}.
	Notice that the authors of \cite{MM13} considered \eqref{two-comp} 
	in the equivalent form with $-t$, $v$ and $u$ instead of $t$, $u$ 
	and $v$, respectively. Moreover, they used the identity
	$\partial_x^2D^{-2}=D^{-2}-I$ in the nonlocal terms, which
	allows to write \eqref{ODE-b} and \eqref{ODE-d} in such a way 
	that the terms with the highest order of derivative involve 
	respectively $\partial_xw$ and $\partial_xz$  only.
	This is crucial for obtaining energy estimates, see \eqref{nonl-term-est} below, 
	but here we left system \eqref{ODE} in the present 
	form to be evidently consistent with \cite[(2.4)]{HM14}.
\end{remark}

\subsection{Solution of the mollified system}
For obtaining the system of ODEs in a Banach space, we consider the following mollification of the problem \eqref{ODE}:
\begin{subequations}
	\label{moll-ODE}
	\begin{align}
	\label{moll-ODE-a}
	&\partial_tu_\ve =-\frac{1}{3}w_\ve^2z_\ve
	+\frac{1}{3}(2u_\ve w_\ve v_\ve+u_\ve^2z_\ve)
	+F(u_\ve,w_\ve,v_\ve,z_\ve),\\
	\nonumber
	&\partial_tw_\ve=
	-J_\ve (w_\ve^\ve(\partial_xw_\ve^\ve)z_\ve^\ve)
	-\frac{1}{3}w_\ve^2v_\ve
	+\frac{4}{3}u_\ve w_\ve z_\ve
	+\frac{1}{3}J_\ve(2u_\ve^\ve(\partial_xw_\ve^\ve)v_\ve^\ve
	+(u_\ve^\ve)^2\partial_xz_\ve^\ve)
	-\frac{2}{3}u_\ve^2v_\ve\\
	\label{moll-ODE-b}
	&\qquad\quad
	+\{J_\ve (u_\ve^\ve(\partial_xw_\ve^\ve)z_\ve^\ve
	-w_\ve^\ve(\partial_xw_\ve^\ve)v_\ve^\ve)
	+u_\ve w_\ve v_\ve
	-u_\ve^2z_\ve\}
	+G_\ve(u_\ve,w_\ve,v_\ve,z_\ve),\\
	\label{moll-ODE-c}
	&\partial_tv_\ve =-\frac{1}{3}w_\ve z_\ve^2
	+\frac{1}{3}(2u_\ve v_\ve z_\ve+w_\ve v_\ve^2)
	+\hat{F}(u_\ve,w_\ve,v_\ve,z_\ve),\\
	\nonumber
	&\partial_tz_\ve=
	-J_\ve (w_\ve^\ve z_\ve^\ve\partial_xz_\ve^\ve)
	-\frac{1}{3}u_\ve z_\ve^2
	+\frac{4}{3}w_\ve v_\ve z_\ve
	+\frac{1}{3}J_\ve(2u_\ve^\ve v_\ve^\ve\partial_xz_\ve^\ve
	+(\partial_xw_\ve^\ve)(v_\ve^\ve)^2)
	-\frac{2}{3}u_\ve v_\ve^2\\
	\label{moll-ODE-d}
	&\qquad\quad
	-\{J_\ve (w_\ve^\ve v_\ve^\ve\partial_xz_\ve^\ve
	-u_\ve^\ve z_\ve^\ve\partial_xz_\ve^\ve)
	+u_\ve v_\ve z_\ve
	-w_\ve v_\ve^2\}
	+\hat{G}_\ve(u_\ve,w_\ve,v_\ve,z_\ve),
	\end{align}
\end{subequations}
with initial data
\begin{equation}
\label{moll-ODE-iv}
u_\ve(x,0)=u_0(x),\quad w_\ve(x,0)=w_0(x),\quad
v_\ve(x,0)=v_0(x),\quad
z_\ve(x,0)=z_0(x),
\end{equation}
and
\begin{subequations}
	\begin{align}
	\nonumber
	& G_\ve=D^{-2}\partial_x\left(
	\frac{1}{3}w_\ve^2z_\ve+
	\left\{u_\ve w_\ve\partial_xz_\ve
	-w_\ve(\partial_xw_\ve)v_\ve
	+\frac{1}{3}J_\ve[u_\ve^\ve
	(u_\ve^\ve\partial_x^2z_\ve^\ve
	-(\partial_x^2w_\ve^\ve)v_\ve^\ve)]
	\right\}
	\right)\\
	\nonumber
	&\qquad +D^{-2}\left(
	\frac{2}{3}u_\ve^2v_\ve+w_\ve^2v_\ve+B(u_\ve,w_\ve,v_\ve,z_\ve)
	\right),\\
	\nonumber
	&\hat G_\ve=D^{-2}\partial_x\left(
	\frac{1}{3}w_\ve z_\ve^2-
	\left\{u_\ve z_\ve\partial_xz_\ve
	-(\partial_xw_\ve)v_\ve z_\ve
	+\frac{1}{3}J_\ve[v_\ve(u_\ve^\ve\partial_x^2z_\ve^\ve
	-(\partial_x^2w_\ve^\ve)v_\ve^\ve)]
	\right\}
	\right)\\
	\nonumber
	&\qquad+D^{-2}\left(
	\frac{2}{3}u_\ve v_\ve^2+u_\ve z_\ve^2
	+\hat B(u_\ve,w_\ve,v_\ve,z_\ve)
	\right).
	\end{align}
\end{subequations}

Notice that in contrast to the FORQ equation \cite{HM14}, $G$ and $\hat G$ 
involve third order derivatives and such terms are mollified 
in $G_\ve$ and $\hat{G}_\ve$.

The system \eqref{moll-ODE} can be interpreted as a Cauchy 
problem in the Banach space 
$$
C\left([-T_\ve,T_\ve],\left(H^{s-1}(A)\right)^4\right),
$$ 
pertaining to the vector $(U_\ve,V_\ve)$, where
\begin{equation}\label{U-V-moll}
U_\ve=(u_\ve,w_\ve),\quad V_\ve=(v_\ve,z_\ve),
\quad
(U_\ve,V_\ve)=(u_\ve,w_\ve,v_\ve,z_\ve),
\end{equation}
with the norm 
\begin{equation}
\nonumber
\|(U_\ve, V_\ve)(t)\|_{H^s}=
\|u_\ve(t)\|_{H^s}+\|w_\ve(t)\|_{H^s}
+\|v_\ve(t)\|_{H^s}+\|z_\ve(t)\|_{H^s}.
\end{equation}
By the Cauchy Theorem in Banach spaces 
(see, e.g., \cite[10.4.5]{D69}) such a problem has a 
unique solution on $[-T_\ve,T_\ve]$ for some $T_\ve>0$.

\begin{proposition}[Energy estimate]
	\label{energy-est}
	Consider the vectors
	$U_\ve(\cdot,t)=(u_\ve(\cdot,t),w_\ve(\cdot,t))$
	and $V_\ve(\cdot,t)=(v_\ve(\cdot,t),z_\ve(\cdot,t))$
	with the norms
	$\|U_\ve(t)\|_{H^s}=\|u_\ve(t)\|_{H^s}
	+\|w_\ve(t)\|_{H^s}$ 
	and
	$\|V_\ve(t)\|_{H^s}=\|v_\ve(t)\|_{H^s}
	+\|z_\ve(t)\|_{H^s}$.
	Then the following inequalities hold, cf.\,\,\cite{HM14}:
	\begin{subequations}\label{energy-U-V}
		\begin{align}
		\label{energy-U}
		&\left|\frac{d}{dt}\|U_\ve(t)\|_{H^{s-1}}\right|
		\lesssim_s
		\|U_\ve(t)\|_{H^{s-1}}^2\|V_\ve(t)\|_{H^{s-1}},
		\quad s>\frac{5}{2},\\
		\label{energy-V}
		&\left|\frac{d}{dt}\|V_\ve(t)\|_{H^{s-1}}\right|
		\lesssim_s
		\|U_\ve(t)\|_{H^{s-1}}\|V_\ve(t)\|_{H^{s-1}}^2,
		\quad s>\frac{5}{2}.
		\end{align}
	\end{subequations}
\end{proposition}

\begin{proof}
	Applying the operator $D^{s-1}$ to \eqref{moll-ODE-b} and 
	then multiplying by $D^{s-1}w_\ve$, on the left hand side we 
	obtain $\frac{1}{2}\frac{d}{dt}\|w_\ve\|_{H^{s-1}}^2$, while
	the first term on the right  hand side can be estimated as 
	follows (recall \eqref{moll-inner-prod} and notice 
	that the operators $D^{s-1}$ and $J_\ve$ commute):
	\begin{align}
	\nonumber
	&\left|\int_A D^{s-1}
	(J_\ve(w_\ve^\ve(\partial_xw_\ve^\ve)z_\ve^\ve))\cdot
	D^{s-1}w_\ve\,dx
	\right|=
	\left|\int_A D^{s-1}
	(w_\ve^\ve(\partial_xw_\ve^\ve)z_\ve^\ve)\cdot
	D^{s-1}w_\ve^\ve\,dx
	\right|\\
	\label{KP-int}
	&\leq\left|
	\int_A\bigl[D^{s-1},w_\ve^\ve z_\ve^\ve\bigr]
	\partial_xw_\ve^\ve\cdot
	D^{s-1}w_\ve^\ve\,dx
	\right|
	+\frac{1}{2}\left|
	\int_A w_\ve^\ve z_\ve^\ve\partial_x(D^{s-1}w_\ve^\ve)^2\,dx
	\right|=:|I_1|+|I_2|,\quad s-1>0.
	\end{align}
	Integrating by parts, using \eqref{Sob-emb}, \eqref{alg} and taking into account \eqref{ineq}, we conclude that $|I_2|$ can be estimated as
	\begin{equation}\label{I2}
	|I_2|\leq\frac{1}{2}\|\partial_x(w_\ve^\ve z_\ve^\ve)\|_{L^\infty}
	\|D^{s-1}w_\ve^\ve\|^2_{L^2}
	\lesssim_s\|\partial_x(w_\ve^\ve z_\ve^\ve)\|_{H^{\hat{s}}}
	\|w_\ve^\ve\|^2_{H^{s-1}}
	\lesssim_s\|w_\ve\|_{H^{\hat{s}+1}}\|z_\ve\|_{H^{\hat{s}+1}}
	\|w_\ve\|^2_{H^{s-1}},
	\end{equation}
	with some $\hat{s}>\frac{1}{2}$. 
	To estimate $|I_1|$ we use, in addition, the Kato-Ponce inequality \eqref{Kato-Ponce}:
	\begin{align}\label{I1}
	\nonumber
	|I_1|&\leq\left\|\left[D^{s-1},w_\ve^\ve z_\ve^\ve\right]
	\partial_xw_\ve^\ve\right\|_{L^2}
	\|D^{s-1}w_\ve^\ve\|_{L^2}\\
	\nonumber
	&\lesssim_s
	(\|D^{s-1}(w_\ve^\ve z_\ve^\ve)\|_{L^2}
	\|\partial_xw_\ve^\ve\|_{L^\infty}
	+\|\partial_x(w_\ve^\ve z_\ve^\ve)\|_{L^\infty}
	\|D^{s-2}\partial_xw_\ve^\ve\|_{L^2})
	\|w_\ve\|_{H^{s-1}}\\
	&\lesssim_s
	\|w_\ve\|^2_{H^{s-1}}\|w_\ve\|_{H^{\hat{s}+1}}\|z_\ve\|_{H^{s-1}}
	+\|w_\ve\|^2_{H^{s-1}}\|w_\ve\|_{H^{\hat{s}+1}}
	\|z_\ve\|_{H^{\hat{s}+1}},
	\end{align}
	with any $\hat{s}>\frac{1}{2}$.
	Taking $\hat{s}=s-2$ in \eqref{I2} and \eqref{I1} we ultimately have 
	\begin{equation}
	\nonumber
	|I_1|+|I_2|\lesssim_s\|w_\ve\|_{H^{s-1}}^3\|z_\ve\|_{H^{s-1}},\quad s>\frac{5}{2}.
	\end{equation}
	The integrals involving 
	$\frac{2}{3}J_\ve(u_\ve^\ve(\partial_xw_\ve^\ve)v_\ve^\ve)$,
	$J_\ve(u_\ve^\ve(\partial_xw_\ve^\ve)z_\ve^\ve)$ and
	$J_\ve(w_\ve^\ve(\partial_xw_\ve^\ve)v_\ve^\ve)$
	can be estimated in a similar manner,
	while for estimating the integral with the nonlocal term
	$D^{-2}\partial_x
	J_\ve(u_\ve^\ve(\partial_x^2w_\ve^\ve)v_\ve^\ve)$ 
	we recall that $D^2=1-\partial_x^2$ and thus
	\begin{equation}\label{nonl-KP}
	D^{-2}\partial_x
	J_\ve(u_\ve^\ve(\partial_x^2w_\ve^\ve)v_\ve^\ve)
	=-J_\ve(u_\ve^\ve(\partial_xw_\ve^\ve)v_\ve^\ve)
	+D^{-2}J_\ve(u_\ve^\ve(\partial_xw_\ve^\ve)v_\ve^\ve)
	-D^{-2}\partial_xJ_\ve(\partial_x(u_\ve^\ve v_\ve^\ve)\partial_xw_\ve^\ve).
	\end{equation}
	The integral with the first term in \eqref{nonl-KP} can be estimated 
	by using the Kato-Ponce inequality, whereas the two other can 
	be treated by the Cauchy-Schwarz inequality and 
	algebra property in $H^{s-2}$ and $H^{s-1}$.
	
	All other terms on the right-hand side of \eqref{moll-ODE-b}, 
	except those involving $\frac{1}{3}J_\ve((u_\ve^\ve)^2\partial_xz_\ve^\ve)$ 
	and $\frac{1}{3}D^{-2}\partial_xJ_\ve
	((u_\ve^\ve)^2(\partial_x^2z_\ve^\ve))$,
	can be estimated by using the Cauchy-Schwarz inequality 
	and the algebra property in $H^{s-2}$ or $H^{s-1}$.
	Thus, eventually we arrive at
	\begin{equation}\label{energy-w-prelim}
	\begin{split}
	&\left|\frac{d}{dt}\|w_\ve\|_{H^{s-1}}^2\right|
	\lesssim_s
	\|w_\ve\|_{H^{s-1}}\|U_\ve\|^2_{H^{s-1}}
	\|V_\ve\|_{H^{s-1}}\\
	&\quad+
	\frac{1}{3}\left|\int_A D^{s-1}
	((u_\ve^\ve)^2\partial_xz_\ve^\ve)\cdot
	D^{s-1}w_\ve^\ve\,dx
	+
	\int_A D^{s-3}\partial_x
	((u_\ve^\ve)^2\partial_x^2z_\ve^\ve)\cdot
	D^{s-1}w_\ve^\ve\,dx\right|,
	\quad s>\frac{5}{2}.
	\end{split}
	\end{equation}
	Recalling that $D^2=1-\partial_x^2$, we conclude
	\begin{equation}\label{nonl-term-est}
	\begin{split}
	\int_A D^{s-3}\partial_x
	((u_\ve^\ve)^2\partial_x^2z_\ve^\ve)\cdot
	D^{s-1}w_\ve^\ve\,dx
	=&
	\int_A D^{s-3}\partial_x^2
	((u_\ve^\ve)^2\partial_xz_\ve^\ve)\cdot
	D^{s-1}w_\ve^\ve\,dx\\
	&-
	\int_A D^{s-3}\partial_x
	(\partial_x(u_\ve^\ve)^2\partial_xz_\ve^\ve)\cdot
	D^{s-1}w_\ve^\ve\,dx\\
	=&
	-\int_A D^{s-1}
	((u_\ve^\ve)^2\partial_xz_\ve^\ve)\cdot
	D^{s-1}w_\ve^\ve\,dx
	+I_3+I_4,
	\end{split}
	\end{equation}
	where
	\begin{align}
	\nonumber
	&I_3=\int_A D^{s-3}
	((u_\ve^\ve)^2\partial_xz_\ve^\ve)\cdot
	D^{s-1}w_\ve^\ve\,dx,\\
	\nonumber
	&I_4=-
	\int_A D^{s-3}\partial_x
	(\partial_x((u_\ve^\ve)^2)\partial_xz_\ve^\ve)\cdot
	D^{s-1}w_\ve^\ve\,dx.
	\end{align}
	Combining \eqref{energy-w-prelim} with \eqref{nonl-term-est}, we 
	conclude that the problematic term $\int D^{s-1}
	((u_\ve^\ve)^2\partial_xz_\ve^\ve)\cdot
	D^{s-1}w_\ve^\ve\,dx$,
	which cannot be treated using the Kato-Ponce inequality together 
	with integration by parts, cf.~\eqref{KP-int}, vanishes. 
	Then, by the Cauchy-Schwarz inequality and the algebra 
	property, we have that
	$|I_3|, |I_4|\lesssim_s\|u_\ve\|_{H^{s-1}}^2\|w_\ve\|_{H^{s-1}}
	\|z_\ve\|_{H^{s-1}}$ for $s>\frac{5}{2}$. 
	Thus we obtain the following estimate
	\begin{equation}\label{energy-w}
	\left|\frac{d}{dt}\|w_\ve\|_{H^{s-1}}\right|
	\lesssim_s
	\|U_\ve\|^2_{H^{s-1}}
	\|V_\ve\|_{H^{s-1}},\quad s>\frac{5}{2}.
	\end{equation}
	Applying the operator $D^{s-1}$ to \eqref{moll-ODE-d}, 
	multiplying by $D^{s-1}z_\ve$ and arguing similarly 
	as for \eqref{moll-ODE-b}, we obtain
	\begin{equation}\label{energy-z}
	\left|\frac{d}{dt}\|z_\ve\|_{H^{s-1}}\right|
	\lesssim_s
	\|U_\ve\|_{H^{s-1}}
	\|V_\ve\|^2_{H^{s-1}},
	\quad s>\frac{5}{2}.
	\end{equation}
	Applying similar techniques as used for equation \eqref{moll-ODE-b}, we 
	can derive estimates for \eqref{moll-ODE-a}, \eqref{moll-ODE-c}:
	\begin{align}
	\label{energy-u}
	&\left|\frac{d}{dt}\|u_\ve\|_{H^{s-1}}\right|
	\lesssim_s
	\|U_\ve\|_{H^{s-1}}^2\|V_\ve\|_{H^{s-1}},
	\quad s>\frac{5}{2},\\
	\label{energy-v}
	&\left|\frac{d}{dt}\|v_\ve\|_{H^{s-1}}\right|
	\lesssim_s
	\|U_\ve\|_{H^{s-1}}\|V_\ve\|_{H^{s-1}}^2,
	\quad s>\frac{5}{2}.
	\end{align}
	Finally, using 
	\eqref{energy-w}, \eqref{energy-u} and 
	\eqref{energy-z}, \eqref{energy-v}, we arrive at 
	\eqref{energy-U} and \eqref{energy-V}, respectively.
\end{proof}

\begin{remark}
	Observe that in the case of the FORQ equation, we have
	$\|U_\ve(t)\|_{H^s}=\|V_\ve(t)\|_{H^s}$,
	and \eqref{energy-U-V} reduces to the energy 
	estimate obtained in \cite{HM14}.
\end{remark}

Using \eqref{energy-U} and \eqref{energy-V} we arrive at 
the following inequality:
\begin{equation}
\nonumber
\left|
\frac{d}{dt}\|(U_\ve,V_\ve)(t)\|_{H^{s-1}}\right|
\leq
C_s\|(U_\ve,V_\ve)(t)\|_{H^{s-1}}^3,\quad 
s>\frac{5}{2},\quad C_s>0,
\end{equation}
which is equivalent to, cf.~\cite[(2.19)]{HM14},
\begin{equation}
\nonumber
\|(U_\ve,V_\ve)(t)\|_{H^{s-1}}\leq
\frac{\|(U_\ve,V_\ve)(0)\|_{H^{s-1}}}
{\sqrt{1-2C_s|t|\|(U_\ve,V_\ve)(0)\|_{H^{s-1}}^2}}.
\end{equation}
Since $\|(U_\ve,V_\ve)(0)\|_{H^{s-1}}\leq
2(\|u_0\|_{H^s}+\|v_0\|_{H^s})$, 
from the latter we have
\begin{equation}
\nonumber
\|(U_\ve,V_\ve)(t)\|_{H^{s-1}}\leq
\frac{2(\|u_0\|_{H^s}+\|v_0\|_{H^s})}
{\sqrt{1-8C_s|t|(\|u_0\|_{H^s}+\|v_0\|_{H^s})^2}}.
\end{equation}
Introducing $T_{\delta_0}=
\frac{1-\delta_0}
{8C_s(\|u_0\|_{H^s}+\|v_0\|_{H^s})^2}$, $0<\delta_0<1$, 
which does not depend on $\ve$, we obtain the following size estimate:
\begin{equation}\label{size-est}
\|(U_\ve,V_\ve)(t)\|_{H^{s-1}}\leq
\frac{2}{\sqrt{\delta_0}}
(\|u_0\|_{H^s}+\|v_0\|_{H^s}),\quad 
-T_{\delta_0}\leq t\leq T_{\delta_0},\quad
0<\delta_0<1.
\end{equation}
Applying \eqref{size-est} together with the algebra property 
of $H^{s-2}$ on the r.h.s.\,\,of \eqref{moll-ODE} we conclude that
\begin{equation}\label{size-est-deriv}
\|\partial_t(U_\ve,V_\ve)(t)\|_{H^{s-2}}
\lesssim_{s,\delta_0}
\|u_0\|_{H^s}^3+\|v_0\|_{H^s}^3,
\quad -T_{\delta_0}\leq t\leq T_{\delta_0},\quad
0<\delta_0<1.
\end{equation}

\subsection{Existence of the solution of (\ref{ODE})}
\label{ex-sol-ODE}
Our objective is to obtain the solution to the Cauchy 
problem \eqref{ODE} using the solution $(U_\ve(t),V_\ve(t))$ 
obtained from the mollified problem \eqref{moll-ODE}-\eqref{moll-ODE-iv}. 
In this subsection, we specifically focus on the case 
where $A=\mathbb{R}$, noting that the periodic problem can 
be handled in a similar manner, as elaborated 
in Remark \ref{ex-periodic} below.

In view of \eqref{size-est}, the family $\{U_\ve,V_\ve\}_{\ve\in(0,1)}$ 
is bounded in $L^\infty\left([-T_{\delta_0},T_{\delta_0}],
\left(H^{s-1}\right)^4\right)$.
Therefore by the Banach-Alaoglu theorem for a sequence $\{\ve_j\}$, $\ve_j\to0$, 
there exists a subsequence $\{\ve_{j_k}\}$ and $(U,V)\in 
L^\infty\left([-T_{\delta_0}, T_{\delta_0}],\left(H^{s-1}\right)^4\right)$ 
such that $\left(U_{\ve_{j_k}},V_{\ve_{j_k}}\right)\to (U,V)$
in the $L^\infty$ weak-$*$ sense, i.e.,
\begin{equation}\label{(U,V)}
\left\langle\left(U_{\ve_{j_k}}, V_{\ve_{j_k}}\right), 
f\right\rangle
\to\langle (U, V), f\rangle,\quad\mbox{for all }
f\in
L^1\left([-T_{\delta_0}, T_{\delta_0}],\left(H^{1-s}\right)^4\right),
\end{equation}
where (cf.\,\,\eqref{U-V-moll})
\begin{equation}
\nonumber
U=(u,w),\quad V=(v,z),
\quad
(U,V)=(u,w,v,z).
\end{equation}
In order to simplify notation, we will write $\ve$ instead of the 
subsequence $\ve_{j_k}$ and set 
$I_{\delta_0}=[-T_{\delta_0}, T_{\delta_0}]$.
Also, we may denote spaces such as 
$L^\infty\left(I_{\delta_0}, (H^s)^4\right)$ as $L^\infty(I_{\delta_0}, H^s)$, 
provided that it does not mislead the reader.

By the $L^\infty$ weak-$*$ lower semicontinuity 
$\|(U,V)\|_{L^\infty(I_{\delta_0},H^{s-1})}
\leq\liminf\limits_{k\to\infty}
\|(U_{\ve_{j_k}}, V_{\ve_{j_k}})\|_{L^\infty(I_{\delta_0},H^{s-1})}$ 
and from \eqref{size-est}, we have the following size 
estimate for the limit $(U,V)$:
\begin{align}
\label{size-est-or}
&\|(U,V)(t)\|_{H^{s-1}}\leq
\frac{2}{\sqrt{\delta_0}}
(\|u_0\|_{H^s}+\|v_0\|_{H^s}),\quad t\in I_{\delta_0}.
\end{align}

In what follows, we will use the smooth cutoff functions
\begin{equation}
\nonumber
\phi_R(x)=
\begin{cases}
1,&x\in[-R+1, R-1],\\
0,&|x|\geq R,
\end{cases}
\end{equation}
where $R>1$ and $0\leq\phi_R(x)\leq1$, which satisfy
the following bound:
\begin{equation}
\nonumber
\max\{
\|\phi_R\|_{H^{s}(\tilde{I}_R)},
\|\phi_R-1\|_{H^{s-1}(\tilde{I}_R)}
\}\leq M,\quad
M>0,\,\,
\tilde{I}_R=(-R,-R+1)\cup(R-1,R),
\end{equation}
with $M$ independent of $R$.
Notice that for such $\phi_R$ we have the estimates
\begin{subequations}\label{phi-est}
	\begin{align}
	\label{phi-est-a}
	&\|\phi_Rf\|_{H^{\hat{s}}(\mathbb{R})}
	\lesssim_{\hat{s}}
	\|f\|_{H^{\hat{s}}([-R+1,R-1])}+
	\|\phi_R\|_{H^{\hat{s}}(\tilde{I}_R)}
	\|f\|_{H^{\hat{s}}(\tilde{I}_R)}
	\lesssim_{\hat{s},M}\|f\|_{H^{\hat{s}}(\mathbb{R})},
	\quad \hat{s}\in(1/2,s-1],\\
	\label{phi-est-b}
	&\|(\partial_x\phi_R)f\|_{H^{\hat{s}}(\mathbb{R})}
	=\|(\partial_x\phi_R)f\|_{H^{\hat{s}}(\tilde{I}_R)}
	\lesssim_{\hat{s},M}\|f\|_{H^{\hat{s}}(\mathbb{R})},
	\quad \hat{s}\in(1/2,s-1],
	\end{align}
\end{subequations}
and
\begin{equation}\label{phi-1-est}
\|(1-\phi_R)f\|_{H^{\hat{s}}(\mathbb{R})}
\lesssim_{\hat{s},M}
\|f\|_{H^{\hat{s}}((-\infty,-R+1]\cup[R-1,\infty))},
\quad \hat{s}\in(1/2,s-1].
\end{equation}

\begin{remark}\label{cutoff-H}
	Observe that if $f,g\in H^{s-1}$ and $\phi_Rf=\phi_Rg$ 
	for all $R>1$, then $f=g$. Indeed, in view of \eqref{phi-1-est}, we 
	can take the limit $R\to\infty$ in the inequality
	$$
	\|f-g\|_{H^{s-1}}\leq
	\|(1-\phi_R)f\|_{H^{s-1}}
	+\|(1-\phi_R)g\|_{H^{s-1}}
	$$
	and obtain that $\|f-g\|_{H^{s-1}}=0$.
\end{remark}

\begin{lemma}\label{comp_UV}
	For any $R>1$ there exists a subsequence 
	of $\phi_R(U_\ve,V_\ve)$ which converges to $\phi_R(U,V)$ in 
	$C\left(I_{\delta_0},H^{s-1-\sigma}\right)$, where $\sigma\in(0,1)$ 
	is such that $s-\sigma>\frac{5}{2}$.
\end{lemma}

\begin{proof} 
	Let us prove that there exists a subsequence $\phi_R U_\ve$ 
	such that $\phi_R U_\ve\to\phi_R U$ in 
	$C\left(I_{\delta_0},H^{s-1-\sigma}\right)$ as $\ve\to 0$.
	We are going to verify the two conditions of the Arzel\`{a}-Ascoli 
	theorem in Banach spaces (see, e.g., \cite[7.5.7]{D69}):
	(i) $\phi_RU_\ve(t)$ is precompact in $H^{s-1-\sigma}$ 
	for all $t\in I_{\delta_0}$ and 
	(ii) $\phi_RU_\ve(t)$ is  equicontinuous in $t$.
	
	Lemma \ref{RK} with $W=[-R,R]$ implies that 
	$H^{s'}([-R,R])\Subset H^{s''}([-R,R])$ for 
	$0<s''<s'<\infty$, $R>0$. Recalling that 
	$U_\ve(t)\in H^{s-1}$ we verify item (i).
	
	From \eqref{phi-est-a} with $\hat{s}=s-1-\sigma$ it follows that 
	$\|\phi_R U_\ve(t_1)-\phi_R U_\ve(t_2)\|_{H^{s-1-\sigma}}
	\lesssim_{s,M}\|U_\ve(t_1)-U_\ve(t_2)\|_{H^{s-1-\sigma}}$ 
	(see \eqref{phi-est}) and therefore 
	for the equicontinuity it is enough to show that
	\begin{equation}\label{U-ve-equicont}
	\|U_\ve(t_1)-U_\ve(t_2)\|_{H^{s-1-\sigma}}
	\lesssim_{s,\delta_0,\|u_0\|_{H^s},\|v_0\|_{H^s}}|t_1-t_2|^{\sigma}.
	\end{equation}
	Taking into account that $a^\sigma\leq 1+a$ for all $a\geq0$ 
	and the energy estimates \eqref{size-est}, \eqref{size-est-deriv}, 
	we conclude that
	\begin{align}
	\nonumber
	\frac{\|U_\ve(t_1)-U_\ve(t_2)\|_{H^{s-1-\sigma}}}
	{|t_1-t_2|^\sigma}&=
	\left(
	\int_{\mathbb{R}}\frac{(1+k^2)^{s-1}}
	{(1+k^2)^\sigma|t_1-t_2|^{2\sigma}}
	(\hat{U}_\ve(t_1)-\hat{U}_\ve(t_2))^2\,dk
	\right)^{1/2}\\
	\nonumber
	&\leq
	\left(
	\int_{\mathbb{R}}
	\left(1+\frac{1}{(1+k^2)(t_1-t_2)^2}\right)
	(1+k^2)^{s-1}
	(\hat{U}_\ve(t_1)-\hat{U}_\ve(t_2))^2\,dk
	\right)^{1/2}\\
	\nonumber
	&\leq\|U_\ve(t_1)-U_\ve(t_2)\|_{H^{s-1}}
	+\frac{\|U_\ve(t_1)-U_\ve(t_2)\|_{H^{s-2}}}{|t_1-t_2|}\\
	\nonumber
	&\leq\frac{2}{\sqrt{\delta_0}}
	(\|u_0\|_{H^s}+\|v_0\|_{H^s})
	+\|\partial_tU_\ve(t_0)\|_{H^{s-2}}\\
	\label{equicont}
	&\lesssim_{s,\delta_0}
	\|u_0\|_{H^s}+\|v_0\|_{H^s}
	+\|u_0\|_{H^s}^3+\|v_0\|_{H^s}^3.
	\end{align}
	and we arrive at \eqref{U-ve-equicont}.
	Thus, we have proved that $\phi_R U_\ve\to\phi_R U$ in
	$C\left(I_{\delta_0},H^{s-1-\sigma}\right)$ as $\ve\to 0$, $\sigma\in(0,1)$.
	Similarly it can be shown that $\phi_R V_\ve\to\phi_R V$ 
	in $C\left(I_{\delta_0},H^{s-1-\sigma}\right)$.
\end{proof}

\begin{corollary}\label{subseq}
	For any $R_1,R_2>1$ there exists a subsequence 
	$(U_{\ve_k},V_{\ve_k})$ such that
	$\phi_{R_j}(U_{\ve_k},V_{\ve_k})$ converges to
	$\phi_{R_j}(U,V)$ in $C\left(I_{\delta_0}, H^{s-1-\sigma}\right)$ 
	as $\ve_k\to0$ for all $j=1,2$.
\end{corollary}

\begin{proof}
	By Lemma \ref{comp_UV} for $R_1$ we have 
	a subsequence $\ve_n$ such that
	$\phi_{R_1}(U_{\ve_n},V_{\ve_n})
	\to\phi_{R_1}(U,V)$ as $\ve_n\to0$.
	Applying Lemma \ref{comp_UV} to the subsequence 
	$(U_{\ve_n},V_{\ve_n})$ and $\phi_{R_2}$ we conclude 
	that there exists $\ve_{n_m}$ such that
	$\phi_{R_2}\left(U_{\ve_{n_m}},V_{\ve_{n_m}}\right)
	\to\phi_{R_2}(U,V)$ as $\ve_{n_m}\to0$.
	Taking $\ve_k=\ve_{n_m}$, we obtain that
	$\phi_{R_j}(U_{\ve_{k}},V_{\ve_{k}})
	\to\phi_{R_j}(U,V)$ as $\ve_k\to0$, $j=1,2$.
	\end{proof}

\begin{proposition}\label{sol-L-infty}
	The pair $(U, V)(t)$ defined by \eqref{(U,V)}
	solves Cauchy problem for the system \eqref{ODE} in 
	$L^{\infty}\left(I_{\delta_0},H^{s-1}\right)$ with 
	initial data $U(0)=(u_0, w_0)$, $V(0)=(v_0, z_0)$.
	
	Moreover, 
	$(U,V)(t)\in\mathrm{Lip}
	\left(I_{\delta_0},H^{s-2}\right)$:
	\begin{equation}\label{Lip-U-V}
	\|(U,V)(t_1)-(U,V)(t_2)\|_{H^{s-2}}
	\lesssim_{s,\delta_0} 
	(\|u_0\|_{H^s}^3+\|v_0\|_{H^s}^3)|t_1-t_2|,
	\quad\mbox{for all }t_1,t_2\in
	I_{\delta_0}.
	\end{equation}
\end{proposition}

\begin{proof}
	We will prove that
	[$\phi_R\,\cdot\,$r.h.s.\,\,of \eqref{moll-ODE}] converges to
	[$\phi_R\,\cdot\,$r.h.s.\,\,of \eqref{ODE}] in 
	$C(I_{\delta_0},H^{s-2-\sigma})$ as $\ve\to0$.
	We give the detailed proof that $\phi_R
	J_\ve(w_\ve^\ve(\partial_xw_\ve^\ve)z_\ve^\ve)$ 
	converges to $\phi_Rw(\partial_xw)z$, the other terms 
	can be treated similarly.
	
	In view of Corollary \ref{subseq} we can take the subsequence 
	$(U_\ve,V_\ve)$ such that $\phi_R(U_\ve,V_\ve)$ and 
	$\phi_{R+1}(U_\ve,V_\ve)$ converge to $\phi_R(U,V)$ 
	and $\phi_{R+1}(U,V)$ respectively.
	Observe that (here $C=C(I_{\delta_0},H^{s-2-\sigma})$)
	\begin{align}\label{moll-to-orig}
	\nonumber
	\|\phi_R
	J_\ve(w_\ve^\ve(\partial_xw_\ve^\ve)z_\ve^\ve)
	-\phi_Rw(\partial_xw)z
	\|_C
	\leq&\,
	\|
	\phi_R
	J_\ve(w_\ve^\ve(\partial_xw_\ve^\ve)z_\ve^\ve)
	-\phi_Rw_\ve^\ve(\partial_xw_\ve^\ve)z_\ve^\ve
	\|_C\\
	&+
	\|\phi_R
	w_\ve^\ve(\partial_xw_\ve^\ve)z_\ve^\ve
	-\phi_Rw(\partial_xw)z
	\|_C.
	\end{align}
	We estimate the first term for all fixed $t$ by using the 
	Sobolev embedding theorem \eqref{Sob-emb}, the algebra 
	property \eqref{alg}, the size estimate \eqref{size-est}, 
	and \eqref{phi-est-a} with $\hat{s}=s-2-\sigma$ as follows:
	\begin{align}
	\nonumber
	\|\phi_R
	J_\ve(w_\ve^\ve(\partial_xw_\ve^\ve)z_\ve^\ve)
	-\phi_R(w_\ve^\ve(\partial_xw_\ve^\ve)z_\ve^\ve)
	\|_{H^{s-2-\sigma}}
	&\lesssim_{s,\sigma,M}
	\|(J_\ve-I)
	(w_\ve^\ve(\partial_xw_\ve^\ve)z_\ve^\ve)\|
	_{H^{s-2-\sigma}}\\
	\nonumber
	&\lesssim_{s,\sigma,M}
	\|(J_\ve-I)\|
	_{\mathcal{L}(H^{s-2},H^{s-2-\sigma})}
	(\|u_0\|_{H^s}^3+\|v_0\|_{H^s}^3),
	\end{align}
	where $\|(J_\ve-I)\|_{\mathcal{L}(H^{s-2},H^{s-2-\sigma})}\to0$ 
	as $\ve\to0$ (see, e.g., \cite[(5.1.35)]{T91}).
	Regarding the second term in \eqref{moll-to-orig} we have
	(here, again, $C=C(I_{\delta_0},H^{s-2-\sigma})$)
	\begin{align}\label{moll-to-orig-2}
	\nonumber
	\|\phi_R
	w_\ve^\ve(\partial_xw_\ve^\ve)z_\ve^\ve
	-\phi_Rw(\partial_xw)z
	\|_C
	\lesssim_{s,\sigma}&\,
	\|w_\ve^\ve z_\ve^\ve\|_C
	\|\phi_R\partial_xw_\ve^\ve
	-\phi_R\partial_xw\|_C\\
	&
	+\|
	\partial_xw
	\|_C
	\|
	\phi_R w_\ve^\ve z_\ve^\ve
	-\phi_R w z
	\|_C.
	\end{align}
	Observe that since $\phi_R\partial_xf=\phi_R\partial_x(\phi_{R+1}f)$, 
	from \eqref{phi-est-a} we have
	\begin{equation}
	\nonumber
	\|\phi_R\partial_xw_\ve^\ve
	-\phi_R\partial_xw\|
	_{C(I_{\delta_0},H^{s-2-\sigma})}
	\lesssim_{s,\sigma,M}
	\|\phi_{R+1}(w_\ve^\ve-w)\|
	_{C(I_{\delta_0},H^{s-1-\sigma})},
	\end{equation}
	and thus 
	$$
	\|\phi_Rw_\ve^\ve v_\ve^\ve-\phi_R wz\|
	_{C}
	\lesssim_{s,\sigma,M}
	\|w_\ve^\ve\|_{C}
	\|\phi_Rz_\ve^\ve-\phi_Rz\|_C+
	\|z\|_C\|\phi_Rw_\ve^\ve-\phi_Rw\|_C,
	$$
	with $C=C(I_{\delta_0},H^{s-2-\sigma})$.
	Using Lemma \ref{comp_UV}, the size estimates 
	\eqref{size-est}, \eqref{size-est-or} and
	$$
	\|\phi_Rf_\ve^\ve-\phi_R f\|_{H^{\hat{s}}}
	\lesssim_{M}\|(J_\ve-I)f_\ve\|_{H^{\hat{s}}}
	+\|\phi_Rf_\ve-\phi_Rf\|_{H^{\hat{s}}},
	\quad\hat{s}\in(1/2,s-1],
	$$
	we conclude that the right-hand side of \eqref{moll-to-orig-2} 
	converges to zero as $\ve\to0$.
	Establishing the convergence of all the other terms
	 by analogy, we can prove that
	$\phi_R\partial_t(U_\ve,V_\ve)
	\to\phi_R\partial_t(U,V)$ in 
	$C(I_{\delta_0},H^{s-2-\sigma})$.
	From Remark \ref{cutoff-H} and that $(U,V)\in L^{\infty}(I_{\delta_0},H^{s-1})$, 
	we conclude that $(U,V)$ is the solution 
	of \eqref{ODE} in $L^{\infty}(I_{\delta_0},H^{s-1})$.
	
	To prove the Lipschitz continuity, observe that 
	in view of \eqref{size-est-deriv} there exists a subsequence (not relabeled) 
	of $\phi_R\partial_t(U_\ve,V_\ve)$ which converges 
	weak-$*$ in $L^\infty(I_{\delta_0},H^{s-2})$:
	\begin{equation}
	\nonumber
	\langle\phi_R
	\partial_t(U_{\ve}, V_{\ve}), f\rangle
	\to\langle
	(\tilde{U}, \tilde{V}), f\rangle,\quad\mbox{for all }
	f\in
	L^1\left(I_{\delta_0},H^{2-s}\right),
	\end{equation}
	for some
	$(\tilde{U},\tilde{V})\in L^\infty(I_{\delta_0},H^{s-2})$.
	Taking into account that $\phi_R\partial_t(U_\ve,V_\ve)\to
	\phi_R\partial_t(U,V)$ weak-$*$ in 
	$L^\infty(I_{\delta_0},H^{s-2-\sigma})$ and that
	$L^1(I_{\delta_0},H^{2-s})\supset
	L^1(I_{\delta_0},H^{2+\sigma-s})$ we conclude that
	$(\tilde{U},\tilde{V})=\phi_R\partial_t(U,V)$. Therefore by the 
	weak-$*$ lower semicontinuity, \eqref{size-est-deriv} 
	and \eqref{phi-est-a} with $\hat{s}=s-2$, the function 
	$\phi_R\partial_t(U,V)\in L^\infty(I_{\delta_0},H^{s-2})$ 
	satisfies the estimate
	\begin{equation}
	\label{size-est-deriv-or}
	\|\phi_R\partial_t(U,V)(t)\|_{H^{s-2}}
	\lesssim_{s,\delta_0,M}
	\|u_0\|_{H^s}^3+\|v_0\|_{H^s}^3,
	\quad t\in I_{\delta_0}.
	\end{equation}
	By the fundamental theorem of calculus in 
	$W^{1,\infty}(I_{\delta_0},H^{s-2})$ we have
	\begin{equation}
	\nonumber
	\phi_R(U,V)(t_1)=\phi_R(U,V)(t_2)
	+\phi_R
	\int_{t_2}^{t_1}\partial_t(U,V)(\tau)\,d\tau,
	\quad\text{for all }t_1,t_2\in I_{\delta_0},
	\end{equation}
	which, combined with \eqref{size-est-deriv-or}, implies the 
	Lipschitz property for $\phi_R(U,V)$:
	\begin{equation}\label{Lip-phi-U-V}
	\|\phi_R(U,V)(t_1)-\phi_R(U,V)(t_2)\|_{H^{s-2}}
	\lesssim_{s,\delta_0,M} 
	(\|u_0\|_{H^s}^3+\|v_0\|_{H^s}^3)|t_1-t_2|,
	\quad t_1,t_2\in
	I_{\delta_0}.
	\end{equation}
	Finally, observe that
	\begin{align}
	\nonumber
	\|(U,V)(t_1)-(U,V)(t_2)\|_{H^{s-2}}\leq&
	\|(1-\phi_R)(U,V)(t_1)\|_{H^{s-2}}
	+\|(1-\phi_R)(U,V)(t_2)\|_{H^{s-2}}\\
	\nonumber
	&+\|\phi_R(U,V)(t_1)-\phi_R(U,V)(t_2)\|_{H^{s-2}}.
	\end{align}
	Using \eqref{phi-1-est} with $\hat{s}=s-2$ and \eqref{Lip-phi-U-V} 
	we conclude that for any $\hat{\ve}>0$ there exists $R>0$ such that
	$$
	\|(U,V)(t_1)-(U,V)(t_2)\|_{H^{s-2}}
	\lesssim_{s,\delta_0,M}
	\hat{\ve}
	+(\|u_0\|_{H^s}^3+\|v_0\|_{H^s}^3)|t_1-t_2|,
	\quad t_1,t_2\in I_{\delta_0},
	$$
	which implies \eqref{Lip-U-V}.
\end{proof}

\begin{remark}
Another way to show convergence of the 
[$\phi_R\,\cdot\,$r.h.s.\,\,of \eqref{moll-ODE}] to
[$\phi_R\,\cdot\,$r.h.s.\,\,of \eqref{ODE}] as $\ve\to0$ is to use the 
Arzel\`{a}-Ascoli theorem in the Banach space
$C(I_{\delta_0},H^{s-2-\sigma})$. 
Here the  Rellich-Kondrachov theorem implies the precompactness condition, while the equicontinuity follows from the size estimate \eqref{size-est} and the estimate \eqref{equicont} (cf.~\cite[Section 4]{HM14}).
\end{remark}

So far we have shown that the solution $(U,V)$ belongs to the space 
$L^{\infty}(I_{\delta_0},H^{s-1})\cap
\mathrm{Lip}(I_{\delta_0},H^{s-2})$, $s>\frac{5}{2}$.
Now our goal is 
to prove that $(U,V)\in C(I_{\delta_0},H^{s-1})$:
\begin{equation}\label{cont-(U,V)}
\lim\limits_{t\to t_0}
\|(U,V)(t)-(U,V)(t_0)\|_{H^{s-1}}^2=0,\quad
\text{for all }t\in I_{\delta_0}.
\end{equation}
However, before proceeding further, we need 
to establish the following lemma:
\begin{lemma}\label{Lip-Y}
	Denote $Y^\ve(t)=
	\|(J_\ve(\phi_R U),
	J_\ve(\phi_R V))(t)\|_{H^{s-1}}^2$.
	Then $Y^\ve(t)$ satisfies the following 
	Lipschitz continuity property on $I_{\delta_0}$:
	\begin{equation}\label{Lip-norm}
	|Y^\ve(t_1)-Y^\ve(t_2)|
	\lesssim
	_{s,\delta_0,\|u_0\|_{H^s},\|v_0\|_{H^s},M}
	|t_1-t_2|,\quad
	t_1,t_2\in I_{\delta_0}.
	\end{equation}
\end{lemma}
\begin{proof}
	By the fundamental theorem of calculus,
	\begin{equation}\label{mvi}
	|Y^\ve(t_1)-Y^\ve(t_2)|\leq
	\sup\limits_{t\in I_{\delta_0}}
	\left|\frac{d}{dt}Y^\ve(t)\right|
	|t_1-t_2|.
	\end{equation}
	After multiplying the equation by $\phi_R$, we apply the 
	mollifier $J_\ve$, followed by the differential 
	operator $D^{s-1}$. Finally, we multiply the 
	result by $D^{s-1}J_\ve(\phi_Ru)$, which gives
	\begin{align}
	\nonumber
	\frac{1}{2}\frac{d}{dt}
	&\|J_\ve(\phi_Ru)\|_{H^{s-1}}^2
	=-\frac{1}{3}\int_{\mathbb{R}} D^{s-1}
	J_\ve(\phi_Rw^2z)\cdot D^{s-1}J_\ve(\phi_Ru)\,dx\\
	\nonumber
	&+\frac{1}{3}\int_{\mathbb{R}} D^{s-1}
	J_\ve(\phi_R(2uwv+u^2z))\cdot 
	D^{s-1}J_\ve(\phi_Ru)\,dx
	+\int_{\mathbb{R}} J_\ve(\phi_RF(u,w,v,z))\cdot 
	D^{s-1}J_\ve(\phi_Ru)\,dx.
	\end{align}
	Using the Cauchy-Schwarz inequality, the algebra property \eqref{alg}, 
	the estimate \eqref{phi-est-a} and that $\|J_\ve f\|_{H^s}\leq\|f\|_{H^s}$, one obtains
	(cf.\,\,\eqref{energy-U-V}; recall that $(U,V)\in L^\infty(I_{\delta_0},H^{s-1})$)
	$$
	\left|\frac{d}{dt}\|J_\ve(\phi_Ru)\|_{H^{s-1}}^2\right|
	\lesssim_{s,M}
	\|U\|_{H^{s-1}}^3\|V\|_{H^{s-1}}.
	$$
	
	The estimate for 
	$\left|\frac{d}{dt}
	\|J_\ve(\phi_R w)\|_{H^{s-1}}^2\right|$ 
	is more involved, due to multipliers 
	like $\partial_x w$ and $\partial_x z$. 
	To deal with the term $w(\partial_x w)z$ we multiply it by $\phi_R$, 
	apply the mollifier $J_\ve$ and the operator $D^{s-1}$, 
	and then multiply the result by $D^{s-1}J_\ve(\phi_R w)$:
	\begin{align}
	\nonumber
	&\left|\int_{\mathbb{R}} D^{s-1}
	J_\ve(\phi_Rw(\partial_xw)z)\cdot
	D^{s-1}J_\ve(\phi_Rw)\,dx
	\right|=
	\left|\int_{\mathbb{R}} D^{s-1}
	(\phi_Rw(\partial_xw)z)\cdot
	D^{s-1}J_\ve^2(\phi_Rw)\,dx
	\right|\\
	\nonumber
	&\leq
	\left|\int_{\mathbb{R}} [D^{s-1},wz]
	\phi_R(\partial_xw)\cdot
	D^{s-1}J_\ve^2(\phi_Rw)\,dx
	\right|
	+
	\left|\int_{\mathbb{R}} wzD^{s-1}
	\phi_R(\partial_xw)\cdot
	D^{s-1}J_\ve^2(\phi_Rw)\,dx
	\right|=I_5+I_6.
	\end{align}
	Here $D^{s-1}
	(\phi_Rw(\partial_xw)z),
	D^{s-1}(\phi_R\partial_xw)\in H^{-1}$, which can be paired 
	with the compactly supported smooth function
	$D^{s-1}J_\ve^2(\phi_Rw)$ (cf.\,\,\cite[Proposition
	5.1.D]{T91};
	notice that for $D^{s-1}J_\ve w$ instead of $D^{s-1}J_\ve (\phi_Rw)$ 
	one can not claim that $D^{s-1}J_\ve^2 w\in H^{1}$).
	
	We can estimate the integral $I_5$ by employing the Kato-Ponce 
	inequality and \eqref{phi-est-a} (cf.~\eqref{I1}):
	\begin{equation}
	\nonumber
	I_5\lesssim_{s,M}\|w\|^3_{H^{s-1}}\|z\|_{H^{s-1}}.
	\end{equation}
	To obtain $\partial_x(D^{s-1}J_\ve(\phi_Rw))^2$ under the 
	integral sign in $I_6$ (cf.~$I_2$ in \eqref{KP-int}) observe that
	\begin{align*}
	I_6\leq&
	\left|
	\int_{\mathbb{R}}J_\ve(wzD^{s-1}\partial_x(\phi_Rw))
	\cdot D^{s-1}J_\ve(\phi_Rw)\,dx
	\right|
	+\left|
	\int_{\mathbb{R}} wz D^{s-1}((\partial_x\phi_R)w)
	\cdot D^{s-1}J_\ve^2(\phi_Rw)\,dx
	\right|\\
	&=I_{61}+I_{62}.
	\end{align*}
	Here, one estimates $I_{62}$ by the Cauchy-Schwarz 
	inequality, Sobolev embedding theorem and \eqref{phi-est-b} 
	with $\hat{s}=s-1$:
	\begin{equation}
	\nonumber
	I_{62}\lesssim_{s,M}
	\|w\|^3_{H^{s-1}}\|z\|_{H^{s-1}},
	\end{equation}
	while $I_{61}$ can be treated as follows:
	\begin{equation}
	\nonumber
	I_{61}\leq
	\left|\int_{\mathbb{R}}
	[J_\ve,wz]D^{s-1}\partial_x(\phi_Rw)
	\cdot D^{s-1}J_\ve(\phi_Rw)\,dx
	\right|
	+\left|\int_{\mathbb{R}}
	wz\partial_x(D^{s-1}J_\ve(\phi_Rw))
	\cdot D^{s-1}J_\ve(\phi_Rw)\,dx
	\right|.
	\end{equation}
	The first term is estimated using the Cauchy-Schwarz 
	inequality and \eqref{HK}, while for the second 
	term we use integration by parts (cf.\,\,\eqref{I2}).
	Thus, one eventually arrives at
	\begin{equation}
	\nonumber
	I_{61}\lesssim_{s,M}
	\|w\|^3_{H^{s-1}}\|z\|_{H^{s-1}}.
	\end{equation}
	The remaining terms on the right-hand side of the equation 
	for $\left|\frac{d}{dt}\|J_\ve \phi_Rw^\ve\|_{H^{s-1}}^2\right|$
	can be handled in a similar manner, where the terms 
	$J_\ve((u_\ve^\ve)^2\partial_xz_\ve^\ve)$ and 
	$D^{-2}\partial_xJ_\ve((u_\ve^\ve)^2(\partial_x^2z_\ve^\ve))$ 
	can be analyzed using \eqref{energy-w-prelim} and \eqref{nonl-term-est}. 
	Eventually, we arrive at the following result:
	\begin{equation}
	\nonumber
	\left|
	\frac{d}{dt}\|J_\ve(\phi_Rw)\|_{H^{s-1}}^2
	\right|
	\lesssim_{s,M}
	\|U\|_{H^{s-1}}^3\|V\|_{H^{s-1}}.
	\end{equation}
	Similarly,
	$\left|
	\frac{d}{dt}\|J_\ve(\phi_Rv)\|_{H^{s-1}}^2
	\right|
	\lesssim_{s,M}
	\! \|U\|_{H^{s-1}}\|V\|_{H^{s-1}}^3$
	and
	$\left|
	\frac{d}{dt}\|J_\ve(\phi_Rz)\|_{H^{s-1}}^2
	\right|
	\lesssim_{s,M}
	\! \|U\|_{H^{s-1}}\|V\|_{H^{s-1}}^3$,
	which, together with the size 
	estimate \eqref{size-est-or} and \eqref{mvi}, imply \eqref{Lip-norm}.
\end{proof}

\begin{proposition}\label{sol-C}
	The solution $(U,V)$ belongs to the space $C(I_{\delta_0},H^{s-1})$.
\end{proposition}

\begin{proof}
	The continuity condition \eqref{cont-(U,V)} of $(U,V)$ is equivalent to
	\begin{align*}\label{cont-Hs-1}
	\lim\limits_{t\to t_0}
	\bigl(&
	\|(U,V)(t)\|_{H^{s-1}}^2
	+\|(U,V)(t_0)\|_{H^{s-1}}^2
	-\langle(U,V)(t),(U,V)(t_0)\rangle_{H^{s-1}}\\
	& \quad -\langle(U,V)(t_0),(U,V)(t)\rangle_{H^{s-1}}
	\bigr)=0,
	\end{align*}
	for all $t_0\in I_{\delta_0}$.
	Using the Lipschitz property \eqref{Lip-U-V} it can be shown that $(U,V)$ 
	is weakly continuous in $H^{s-1}$ (see, e.g., \cite[Lemma 2]{HM14}).
	Therefore it is enough to show that
	$\lim\limits_{t\to t_0}\|(U,V)(t)\|_{H^{s-1}}^2
	=\|(U,V)(t_0)\|_{H^{s-1}}^2$.
	Notice that from Lemma \ref{Lip-Y}, \eqref{phi-1-est} and the inequality
	\begin{equation}
	\nonumber
	\begin{split}
	\|f^\ve(t_1)-f^\ve(t_2)\|^2_{H^{s-1}}
	\lesssim&\,
	\|J_\ve(\phi_R-1)f(t_1)\|^2_{H^{s-1}}
	+\|J_\ve(\phi_R-1)f(t_2)\|^2_{H^{s-1}}\\
	&+\|J_\ve(\phi_Rf)(t_1)
	-J_\ve(\phi_Rf)(t_2)\|^2_{H^{s-1}},
	\end{split}
	\end{equation}
	we conclude that $\|(U^\ve,V^\ve)(t)\|_{H^{s-1}}^2$ 
	is Lipschitz continuous on $I_{\delta_0}$.
	Taking into account that 
	$$
	\|(U^\ve,V^\ve)(t)\|_{H^{s-1}}\to\|(U,V)(t)\|_{H^{s-1}},
	\quad \text{as $\ve\to0$},
	$$ 
	for all fixed $t\in I_{\delta_0}$
	we conclude from the Lipschitz property of the 
	map $t\mapsto\|(U^\ve,V^\ve)(t)\|_{H^{s-1}}^2$ that 
	the map $t\mapsto\|(U,V)(t)\|_{H^{s-1}}^2$ also 
	satisfies the Lipschitz property and, consequently, both 
	maps $t\mapsto\|(U,V)(t)\|_{H^{s-1}}^2$ 
	and $t\mapsto\|(U,V)(t)\|_{H^{s-1}}$ are continuous.
\end{proof}

\begin{remark}\label{ex-periodic}
	Proving existence of solutions for the periodic 
	problem (i.e., the case $A=\mathbb{R}/2\pi\mathbb{Z}$) requires 
	essentially one modification: there is no need to use the cutoff function $\phi_R$.
	Indeed, since the circle is a compact set, we have in the proof 
	of Lemma \ref{comp_UV} that $H^{s-1-\sigma}(A)\Subset H^{s-1}(A)$ 
	directly from the Rellich-Kondrachov theorem.
	Thus in Lemma \ref{comp_UV} one obtains that a subsequence 
	of $(U_\ve,V_\ve)$ converges to $(U,V)$, in Lemma \ref{Lip-Y} one 
	should take $Y^\ve(t)=\|(U^\ve,V^\ve)(t)\|_{H^{s-1}}^2$, while
	the proofs of Propositions \ref{sol-L-infty} and \ref{sol-C} become 
	simpler in view of the absence of the factor $\phi_R$.
\end{remark}

Combining Propositions \ref{sol-L-infty}, \ref{sol-C} and 
Remark \ref{ex-periodic}, we conclude that for 
initial data $u_0=u(0)$, $v_0=v(0)$ belonging 
to $H^{s}(A)$, $s>\frac{5}{2}$, there exists a solution 
$(U,V)\in C([-T_{\delta_0},T_{\delta_0}],H^{s-1}(A))$ 
of the Cauchy problem for \eqref{ODE}, where
$T_{\delta_0}=\frac{1-\delta_0}
{8C_s(\|u_0\|_{H^s}+\|v_0\|_{H^s})^2}$, with 
some $C_s>0$ and $0<\delta_0<1$.

\subsection{Uniqueness of the solution of (\ref{ODE})}
In this subsection we prove that the solution of \eqref{ODE} is unique.
Assume that there exist two solutions $(u,w,v,z)$, 
$(\tilde{u},\tilde{w},\tilde{v},\tilde{z})$ and then set
\begin{equation}\label{phi-psi}
\phi_1=u-\tilde{u},\,\,
\psi_1=w-\tilde{w},\,\,
\phi_2=v-\tilde{v},\,\,
\psi_2=z-\tilde{z}.
\end{equation}
Subtracting the equations for $(u,w,v,z)$ and 
$(\tilde{u},\tilde{w},\tilde{v},\tilde{z})$, using the identities
$$f^2g-\tilde{f}^2\tilde{g}=f^2(g-\tilde{g})
+\tilde{g}(f+\tilde{f})(f-\tilde{f}),\quad
fgh-\tilde{f}\tilde{g}\tilde{h}=
fg(h-\tilde{h})+f\tilde{h}(g-\tilde{g})
+\tilde{g}\tilde{h}(f-\tilde{f}),$$
we obtain the following differential equations for $\phi_j$, $\psi_j$, $j=1,2$:
\begin{subequations}
	\label{ODE-phi-psi}
	\begin{align}
	\label{ODE-phi-psi-a}
	&\partial_t\phi_1=
	\alpha_1\phi_1
	+\alpha_2\psi_1
	+\frac{2}{3}uw\phi_2
	+\frac{1}{3}(u^2-w^2)\psi_2
	+F(u,w,v,z)
	-F(\tilde{u},\tilde{w},\tilde{v},\tilde{z}),\\
	\label{ODE-phi-psi-b}
	&\partial_t\psi_1=
	\beta_1\phi_1
	+\beta_2\partial_x\psi_1
	+\beta_3\psi_1
	+\beta_4\phi_2
	+\frac{1}{3}u^2\partial_x\psi_2
	+\beta_5\psi_2
	+G(u,w,v,z)
	-G(\tilde{u},\tilde{w},\tilde{v},\tilde{z}),\\
	&\partial_t\phi_2=
	\frac{2}{3}\tilde{v}\tilde{z}\phi_1
	+\frac{1}{3}(\tilde{v}^2-\tilde{z}^2)\psi_1
	+\gamma_1\phi_2
	+\gamma_2\psi_2
	+\hat F(u,w,v,z)
	-\hat F(\tilde{u},\tilde{w},\tilde{v},\tilde{z}),\\
	&\partial_t\psi_2=
	\delta_1\phi_1
	+\frac{1}{3}v^2\partial_x\psi_1
	+\delta_2\psi_1
	+\delta_3\phi_2
	+\delta_4\partial_x\psi_2
	+\delta_5\psi_2
	+\hat G(u,w,v,z)
	-\hat G(\tilde{u},\tilde{w},\tilde{v},\tilde{z}),
	\end{align}
\end{subequations}
where
\begin{equation}
\label{coeff-loc}
\begin{split}
&\alpha_1=\frac{1}{3}(u+\tilde{u})\tilde{z}
+\frac{2}{3}\tilde{w}\tilde{v},\quad
\alpha_2=\frac{2}{3}u\tilde{v}
-\frac{1}{3}(w+\tilde{w})\tilde{z},\\
&\beta_1=(u+\tilde{u})
\left(
-\frac{2}{3}\tilde{v}
+\frac{1}{3}\partial_x\tilde{z}
-\tilde{z}
\right)
+(\partial_x\tilde{w})
\left(
\frac{2}{3}\tilde{v}+\tilde{z}\right)
+\tilde{w}\left(\tilde{v}+\frac{4}{3}\tilde{z}
\right),\\
&\beta_2=(u-w)\tilde{z}
+\tilde{v}\left(\frac{2}{3}u-w\right),\quad
\beta_3=u\left(\tilde{v}+\frac{4}{3}\tilde{z}\right)
-(\partial_x\tilde{w})(\tilde{v}+\tilde{z})
-\frac{1}{3}(w+\tilde{w})\tilde{v},\\
&\beta_4=\frac{2}{3}u\left(
-u+\partial_xw+\frac{3}{2}w
\right)
-w\left(\partial_xw+\frac{1}{3}w\right),\quad
\beta_5=u\left(
-u+\partial_xw+\frac{4}{3}w\right)
-w\partial_xw,\\
&\gamma_1=\frac{1}{3}(2u\tilde{z}
+\tilde{w}(v+\tilde{v})),\quad
\gamma_2=
\frac{2}{3}uv-\frac{1}{3}\tilde{w}(z+\tilde{z}),\\
&\delta_1=-\frac{2}{3}v^2
+\tilde{v}\left(
\frac{2}{3}\partial_x\tilde{z}-\tilde{z}
\right)+\tilde{z}\partial_x\tilde{z}
-\frac{1}{3}z^2,\quad
\delta_2=
v^2-\tilde{v}\left(\partial_x\tilde{z}-
\frac{4}{3}\tilde{z}\right)
-\tilde{z}\partial_x\tilde{z},\\
&\delta_3=
u\left(
\frac{2}{3}\partial_x\tilde{z}-\tilde{z}
\right)
-(v+\tilde{v})\left(
\frac{2}{3}\tilde{u}
-\frac{1}{3}\partial_x\tilde{w}-\tilde{w}
\right)
-w\left(
\partial_x\tilde{z}-\frac{4}{3}\tilde{z}
\right),\\
&\delta_4=\left(
\frac{2}{3}u-w
\right)v
+\left(
u-w
\right)z,\quad
\delta_5=
-u\left(
v-\partial_x\tilde{z}
\right)
-\frac{1}{3}\tilde{u}(z+\tilde{z})
+w\left(
\frac{4}{3}v-\partial_x\tilde{z}
\right).
\end{split}
\end{equation}
The differences of the nonlocal terms can be written as follows:
\begin{subequations}
	\begin{align}
	\nonumber
	&F(u,w,v,z)
	-F(\tilde{u},\tilde{w},\tilde{v},\tilde{z})=
	D^{-2}\left(
	\eta_1\phi_1-\frac{1}{3}u\tilde{v}\partial_x^2\psi_1
	-w\tilde{v}\partial_x\psi_1+\eta_2\psi_1+\eta_3\phi_2
	+\frac{1}{3}u^2\partial_x^2\psi_2\right.\\
	\label{F-phi-psi}
	&\left.\quad+uw\partial_x\psi_2
	+\frac{1}{3}w^2\psi_2
	\right)
	+D^{-2}\partial_x
	\bigg(
	\eta_4\phi_1
	+\eta_5\psi_1
	+\eta_6\phi_2
	+B(u,w,v,z)
	-B(\tilde{u},\tilde{w},\tilde{v},\tilde{z})
	\bigg),\\
	\nonumber
	&G(u,w,v,z)
	-G(\tilde{u},\tilde{w},\tilde{v},\tilde{z})=
	D^{-2}\partial_x
	\left(
	\eta_1\phi_1-\frac{1}{3}u\tilde{v}\partial_x^2\psi_1
	-w\tilde{v}\partial_x\psi_1+\eta_2\psi_1+\eta_3\phi_2
	+\frac{1}{3}u^2\partial_x^2\psi_2\right.\\
	\label{G-phi-psi}
	&\left.\quad+uw\partial_x\psi_2
	+\frac{1}{3}w^2\psi_2
	\right)
	+D^{-2}\bigg(
	\eta_4\phi_1
	+\eta_5\psi_1
	+\eta_6\phi_2
	+B(u,w,v,z)
	-B(\tilde{u},\tilde{w},\tilde{v},\tilde{z})
	\bigg),
	\end{align}
\end{subequations}
	and
	\begin{align}
	\nonumber
	&\hat F(u,w,v,z)
	-\hat F(\tilde{u},\tilde{w},\tilde{v},\tilde{z})=
	D^{-2}\left(
	\nu_1\phi_1+\frac{1}{3}v^2\partial_x^2\psi_1
	-\tilde{v}\tilde{z}\partial_x\psi_1
	+\frac{1}{3}z^2\psi_1+\nu_2\phi_2
	-\frac{1}{3}uv\partial_x^2\psi_2\right.\\
	\nonumber
	&\left.\quad-uz\partial_x\psi_2
	+\nu_3\psi_2
	\right)
	+D^{-2}\partial_x\left(
	\nu_4\phi_1
	+\nu_5\phi_2
	+\nu_6\psi_2
	+\hat{B}(u,w,v,z)
	-\hat{B}(\tilde{u},\tilde{w},\tilde{v},\tilde{z})
	\right),\\
	\nonumber
	&\hat G(u,w,v,z)
	-\hat G(\tilde{u},\tilde{w},\tilde{v},\tilde{z})=
	D^{-2}\partial_x\left(
	\nu_1\phi_1+\frac{1}{3}v^2\partial_x^2\psi_1
	-\tilde{v}\tilde{z}\partial_x\psi_1
	+\frac{1}{3}z^2\psi_1+\nu_2\phi_2
	-\frac{1}{3}uv\partial_x^2\psi_2\right.\\
	\nonumber
	&\left.\quad-uz\partial_x\psi_2
	+\nu_3\psi_2
	\right)
	+D^{-2}\left(
	\nu_4\phi_1
	+\nu_5\phi_2
	+\nu_6\psi_2
	+\hat{B}(u,w,v,z)
	-\hat{B}(\tilde{u},\tilde{w},\tilde{v},\tilde{z})
	\right),
	\end{align}
with
	\begin{align*}
	&B(u,w,v,z)
	-B(\tilde{u},\tilde{w},\tilde{v},\tilde{z})=
	\eta_7\phi_1
	+\eta_8\partial_x\psi_1
	+\eta_9\psi_1
	+\eta_{10}\phi_2
	-\frac{1}{3}w^2\partial_x\psi_2
	+\eta_{11}\psi_2,
	\\
	&\hat{B}(u,w,v,z)
	-\hat{B}(\tilde{u},\tilde{w},\tilde{v},\tilde{z})=
	\nu_7\phi_1
	-\frac{1}{3}z^2\partial_x\psi_1
	+\nu_8\psi_1
	+\nu_9\phi_2
	+\nu_{10}\partial_x\psi_2
	+\nu_{11}\psi_2,
	\end{align*}
where
\begin{equation}\label{eta}
\begin{split}
&\eta_1=
\frac{1}{3}(u+\tilde{u})\partial_x^2\tilde{z}
-\frac{1}{3}(\partial_x^2\tilde{w})\tilde{v}
+\tilde{w}\partial_x\tilde{z},\quad
\eta_2=u\partial_x\tilde{z}
-(\partial_x\tilde{w})\tilde{v}
+\frac{1}{3}(w+\tilde{w})\tilde{z},\\
&\eta_3=-\frac{1}{3}u\partial_x^2w
-w\partial_xw,\quad
\eta_4=\frac{2}{3}\left(
u+\tilde{u}
\right)\tilde{v},\quad
\eta_5=(w+\tilde{w})\tilde{v},\quad
\eta_6=\frac{2}{3}u^2+w^2,\\
&\eta_7=(u+\tilde{u})\tilde{z}
-(\partial_x\tilde{w})\tilde{z}
-\tilde{w}\tilde{v},\quad
\eta_8=-u\tilde{z}
+w\left(\tilde{v}+\frac{1}{3}\tilde{z}\right),\\
&\eta_9=-u\tilde{v}
+(\partial_x\tilde{w})
\left(\tilde{v}+\frac{1}{3}\tilde{z}\right)
-\frac{1}{3}(\partial_x\tilde{z})(w+\tilde{w}),\quad
\eta_{10}=-uw+w\partial_xw,\\
&\eta_{11}=u(u-\partial_xw)+\frac{1}{3}w\partial_xw,
\end{split}
\end{equation}
and
\begin{equation}
\nonumber
\begin{split}
&\nu_1=-\frac{1}{3}\tilde{v}\partial_x^2\tilde{z}
-\tilde{z}\partial_x\tilde{z},\quad
\nu_2=-\frac{1}{3}u\partial_x^2\tilde{z}
+\frac{1}{3}(\partial_x^2\tilde{w})(v+\tilde{v})
+(\partial_xw)\tilde{z},\\
&\nu_3=-u\partial_x\tilde{z}+(\partial_xw)v
+\frac{1}{3}\tilde{w}(z+\tilde{z}),\quad
\nu_4=\frac{2}{3}v^2+z^2,\quad
\nu_5=\frac{2}{3}\tilde{u}(v+\tilde{v}),\quad
\nu_6=\tilde{u}(z+\tilde{z}),\\
&\nu_7=\tilde{v}\tilde{z}-\tilde{z}\partial_x\tilde{z},
\quad
\nu_8=-v^2+\tilde{v}\partial_x\tilde{z}
+\frac{1}{3}\tilde{z}\partial_x\tilde{z},\quad
\nu_9=u\tilde{z}-\tilde{w}(v+\tilde{v})
+w\partial_x\tilde{z},\\
&\nu_{10}=-uz+wv+\frac{1}{3}wz,\quad
\nu_{11}=u(v-\partial_x\tilde{z})
-\frac{1}{3}\partial_x\tilde{w}(z+\tilde{z})
+\frac{1}{3}w\partial_x\tilde{z}.
\end{split}
\end{equation}
\begin{proposition}\label{energy-phi-psi}
	Let $\sigma\in(1/2,s-2]$.
	Then we have the following energy estimates 
	for $\phi_j$ and $\psi_j$, $j=1,2$:
	\begin{subequations}
		\begin{align}
		\label{energy-phi}
		&\left|\frac{d}{dt}\|\phi_j(t)\|_{H^{\sigma}}\right|
		\lesssim_{\sigma,s,\delta_0}
		\left(\|u_0\|^2_{H^s}+\|v_0\|^2_{H^s}\right)
		\sum\limits_{i=1}^2\bigl(
		\|\phi_i(t)\|_{H^{\sigma}}
		+\|\psi_i(t)\|_{H^{\sigma}}\bigr),
		\quad t\in I_{\delta_0},\,\,j=1,2,\\
		\label{energy-psi}
		&\left|\frac{d}{dt}\|\psi_j(t)\|_{H^{\sigma}}\right|
		\lesssim_{\sigma,s,\delta_0}
		\left(\|u_0\|^2_{H^s}+\|v_0\|^2_{H^s}\right)
		\sum\limits_{i=1}^2\bigl(
		\|\phi_i(t)\|_{H^{\sigma}}
		+\|\psi_i(t)\|_{H^{\sigma}}\bigr),
		\quad t\in I_{\delta_0},\,\,j=1,2.
		\end{align}
	\end{subequations}
\end{proposition}

\begin{proof}
	First, let us prove the energy estimate for $\phi_1$.
	Applying the operator $D^\sigma$ to \eqref{ODE-phi-psi-a} 
	and multiplying by $D^\sigma\phi_1$, the left-hand side becomes 
	$\frac{1}{2}\frac{d}{dt}\|\phi_1\|_{H^\sigma}^2$.
	The modulus of all the integrals on the right-hand side, except the ones involving
	$\partial_x^2w$, $\partial_x^2\tilde{w}$, $\partial_x^2\tilde{z}$ in 
	the nonlocal terms (see $\eta_1$ and $\eta_3$ in \eqref{eta}) as well as
	$\partial_x\psi_j$ and $\partial_x^2\psi_j$, $j=1,2$ (see \eqref{F-phi-psi}), 
	can be estimated using the Cauchy-Schwarz inequality, the 
	algebra property in either $H^{\sigma}$ or $H^{\sigma+1}$ and 
	the size estimate \eqref{size-est-or}, where we take into 
	account that $\sigma+1\leq s-1$.
	
	Let us show how to estimate the integral which involves 
	$\frac{1}{3}u\partial_x^2\tilde{z}$,
	see $\eta_1$ in \eqref{eta}; the terms 
	$\frac{1}{3}\tilde{u}\partial_x^2\tilde{z}$, 
	$\frac{1}{3}(\partial_x^2\tilde{w})\tilde{v}$ and
	$\frac{1}{3}u\partial_x^2w$, see respectively $\eta_1$ and $\eta_3$, 
	can be treated similarly. 
	By the Cauchy-Schwarz inequality we have
	\begin{equation*}
	\left|
	\int_A D^{\sigma-2}
	(u(\partial_x^2\tilde{z})\phi_1)
	\cdot D^\sigma\phi_1\,dx
	\right|\leq
	\|(u\partial_x^2\tilde{z}\phi_1)\|
	_{H^{\sigma-2}}
	\|\phi_1\|_{H^{\sigma}}.
	\end{equation*}
	Taking into account that (recall $\partial_x\phi_1=\psi_1$, 
	see Remark \ref{u_x=w})
	\begin{equation}
	\nonumber
	u(\partial_x^2\tilde{z})\phi_1=
	\partial_x(u(\partial_x\tilde{z})\phi_1)
	-(\partial_xu)(\partial_x\tilde{z})\phi_1
	-u(\partial_x\tilde{z})\psi_1,
	\end{equation}
	we obtain
	\begin{equation}\label{two-deriv-z}
	\|(u(\partial_x^2\tilde{z})\phi_1)\|
	_{H^{\sigma-2}}
	\leq
	\|u(\partial_x\tilde{z})\phi_1\|
	_{H^{\sigma}}
	+\|(\partial_xu)
	(\partial_x\tilde{z})\phi_1\|_{H^{\sigma}}
	+\|u(\partial_x\tilde{z})\psi_1\|
	_{H^{\sigma}}.
	\end{equation}
	The norms on the right-hand side can be easily estimated by the 
	algebra property and the size estimate \eqref{size-est-or}.
	
	Let us estimate the integral which involves 
	$D^{-2}\partial_x(w\tilde{z}\partial_x\psi_1)$, see $\eta_8$ in \eqref{eta}
	(all other integrals which involve terms in 
	$D^{-2}\partial_x(\eta_8\partial_x\psi_1)$ 
	as well as $D^{-2}\partial_x(w^2\partial_x\psi_2)$,
	$D^{-2}(w\tilde{v}\partial_x\psi_1)$ and
	$D^{-2}(uw\partial_x\psi_2)$
	can be estimated similarly).
	Using the Cauchy-Schwarz inequality and that 
	$\|\partial_xf\|_{H^{\sigma-2}}\leq
	\|f\|_{H^{\sigma-1}}$, we obtain
	\begin{equation*}
	\left|
	\int_A D^{\sigma-2}\partial_x
	(w\tilde{z}\partial_x\psi_1)
	\cdot D^\sigma\phi_1\,dx
	\right|\leq
	\|w\tilde{z}\partial_x\psi_1\|
	_{H^{\sigma-1}}
	\|\phi_1\|_{H^{\sigma}}.
	\end{equation*}
	Taking into account that $\|f\partial_x\psi_1\|_{H^{\sigma-1}}
	\leq\|f\psi_1\|_{H^{\sigma}}
	+\|(\partial_xf)\psi_1\|_{H^{\sigma-1}}$ together with the algebra 
	property in $H^\sigma$ and the size estimate (recall 
	that $\sigma+1\leq s-1$), we obtain the desired inequality.
	
	Finally, let us show how to deal with 
	$\frac{1}{3}u\tilde{v}\partial_x^2\psi_1$, the term 
	$\frac{1}{3}u^2\partial_x^2\psi_2$ can be handled similarly.
	By the Cauchy-Schwarz inequality, 
	\begin{equation*}
	\left|
	\int_A D^{\sigma-2}(u\tilde{v}\partial_x^2\psi_1)
	\cdot D^\sigma\phi_1\,dx
	\right|\leq
	\|(u\tilde{v}\partial_x^2\psi_1)\|
	_{H^{\sigma-2}}
	\|\phi_1\|_{H^{\sigma}}.
	\end{equation*}
	Taking into account that
	\begin{equation*}
	u\tilde{v}\partial_x^2\psi_1
	=\partial_x^2(u\tilde{v}\psi_1)
	-2\partial_x[(\partial_x(u\tilde{v}))\psi_1]
	+(\partial_x^2(u\tilde{v}))\psi_1,
	\end{equation*}
	and using $\partial_xu=w$ and 
	$\partial_x\tilde{v}=\tilde{z}$ (see Remark \ref{u_x=w}) we obtain
	\begin{equation}\label{two-deriv-psi}
	\|(u\tilde{v}\partial_x^2\psi_1)\|
	_{H^{\sigma-2}}
	\lesssim
	\|u\tilde{v}\psi_1\|_{H^{\sigma}}
	+\|(\partial_x(u\tilde{v}))\psi_1\|_{H^\sigma}
	+\|(\partial_x(w\tilde{v}+u\tilde{z}))\psi_1\|
	_{H^\sigma},
	\end{equation}
	which can be further estimated by the algebra property and 
	the size estimate \eqref{size-est-or}.
	Considering the energy estimate for $\phi_2$ similarly, we 
	arrive at \eqref{energy-phi}.
	
	Now let us prove \eqref{energy-psi} for $\psi_1$.
	We apply the operator $D^\sigma$ to \eqref{ODE-phi-psi-b} 
	and multiply by $D^\sigma\psi_1$.
	On the left-hand side we obtain 
	$\frac{1}{2}\frac{d}{dt}\|\psi_1\|^2_{H^\sigma}$.
	All the integral terms on the right-hand side can be estimated 
	similarly to $\phi_1$, except for those involving
	$\beta_2\partial_x\psi_1$, 
	$\frac{1}{3}u^2\partial_x\psi_2$ as well as
	$\frac{1}{3}D^{-2}
	\partial_x(u\tilde{v}\partial_x^2\psi_1)$
	and
	$\frac{1}{3}D^{-2}
	\partial_x(u^2\partial_x^2\psi_2)$ (see \eqref{G-phi-psi}).
	Let us show how to estimate the integral with the term 
	$w\tilde{z}\partial_x\psi_1$ which is present in $\beta_2\partial_x\psi_1$, 
	see \eqref{coeff-loc}:
	\begin{equation}\label{CCM-est}
	\begin{split}
	\left|
	\int_A D^\sigma(w\tilde{z}\partial_x\psi_1)\cdot
	D^\sigma\psi_1\,dx	\right|
	\leq&
	\left|
	\int_A [D^\sigma\partial_x,w\tilde{z}]\psi_1\cdot
	D^\sigma\psi_1\,dx
	\right|
	+
	\left|
	\int_A w\tilde{z}D^\sigma\partial_x\psi_1\cdot
	D^\sigma\psi_1\,dx
	\right|\\
	&+
	\left|
	\int_A D^\sigma[(\partial_x(w\tilde{z}))\psi_1]\cdot
	D^\sigma\psi_1\,dx
	\right|
	=|I_7|+|I_8|+|I_9|.
	\end{split}
	\end{equation}
	Here, $|I_9|$ can be estimated by the Cauchy-Schwarz 
	inequality, the algebra property in $H^\sigma$ and the size 
	estimate \eqref{size-est-or}, where we use that $\sigma+1\leq s-1$.
	We estimate $|I_8|$ using integration by parts:
	\begin{equation}
	\nonumber
	\begin{split}
	|I_8|&=
	\frac{1}{2}\left|
	\int_A w\tilde{z}\partial_x(D^\sigma\psi_1)^2\,dx
	\right|
	=
	\frac{1}{2}\left|
	\int_A (\partial_x(w\tilde{z}))(D^\sigma\psi_1)^2\,dx
	\right|
	\lesssim\|\partial_x(w\tilde{z})\|_{L^2}
	\|\psi_1\|^2_{H^\sigma}\\
	&\lesssim
	\|w\tilde{z}\|_{H^{s-1}}
	\|\psi_1\|^2_{H^\sigma},
	\end{split}
	\end{equation}
	and using the size estimate \eqref{size-est-or}.
	To deal with $|I_7|$ we use the Cauchy-Schwarz inequality 
	and the Calderon-Coifman-Meyer estimate \eqref{CCM}, 
	which gives us the following inequality
	\begin{equation}\label{I_7}
	|I_7|\leq
	\|[D^\sigma\partial_x,w\tilde{z}]\psi_1\|_{L^2}
	\|\psi_1\|_{H^\sigma}
	\lesssim_{s,\sigma}
	\|w\|_{H^{s-1}}\|\psi_1\|^2_{H^\sigma}.
	\end{equation}
	All the other integrals which involve the terms in $\beta_2\partial_x\psi_1$
	can be estimated similarly.
	Concerning the nonlocal term $\frac{1}{3}D^{-2}
	\partial_x(u\tilde{v}\partial_x^2\psi_1)$ we 
	observe, cf.~\eqref{nonl-KP} (recall that $\partial_xu=w$ and 
	$\partial_x\tilde{v}=\tilde{z}$),
	\begin{equation}
	\nonumber
	D^{-2}
	\partial_x(u\tilde{v}\partial_x^2\psi_1)=
	-u\tilde{v}\partial_x\psi_1
	+D^{-2}(u\tilde{v}\partial_x\psi_1)
	-D^{-2}\partial_x(w\tilde{v}\partial_x\psi_1)
	-D^{-2}\partial_x(u\tilde{z}\partial_x\psi_1),
	\end{equation}
	which implies that the first term can be estimated as 
	in \eqref{CCM-est}, while the other terms can be 
	treated by the Cauchy-Schwarz inequality and
	\begin{equation}\label{pd-psi}
	\begin{split}
	\|w\tilde{v}\partial_x\psi_1\|_{H^{\sigma-1}}
	&\leq\|w\tilde{v}\psi_1\|_{H^{\sigma}}
	+\|(\partial_xw)\tilde{v}\psi_1\|_{H^{\sigma-1}}
	+\|w\partial_x\tilde{v}\psi_1\|_{H^{\sigma-1}}\\
	&\lesssim_\sigma
	(\|w\tilde{v}\|_{H^{\sigma}}
	+\|(\partial_xw)\tilde{v}\|_{H^{\sigma}}
	+\|w\partial_x\tilde{v}\|_{H^{\sigma}})
	\|\psi_1\|_{H^{\sigma}}
	\lesssim_{\sigma,s}
	\|w\tilde{v}\|_{H^{s-1}}\|\psi_1\|_{H^{\sigma}}.
	\end{split}
	\end{equation}
	
	It remains to estimate the integrals involving 
	$\frac{1}{3}u^2\partial_x\psi_2$ and $\frac{1}{3}D^{-2}
	\partial_x(u^2\partial_x^2\psi_2)$ (notice that the former cannot 
	be treated by following 
	\eqref{CCM-est}--\eqref{I_7}, since we are not able to 
	integrate by parts). Arguing similarly as in \eqref{nonl-term-est}, 
	we conclude that
	\begin{equation}
	\nonumber
	\begin{split}
	&\left|
	\int_A D^\sigma(u^2\partial_x\psi_2)
	\cdot D^\sigma\psi_1\,dx
	+\int_A D^{\sigma-2}\partial_x(u^2\partial_x^2\psi_2)
	\cdot D^\sigma\psi_1\,dx
	\right|\\
	&\qquad\leq\left|
	\int_A D^{\sigma-2}(u^2\partial_x\psi_2)
	\cdot D^\sigma\psi_1\,dx
	\right|
	+2\left|
	\int_A D^{\sigma-2}
	\partial_x(u(\partial_xu)\partial_x\psi_2)
	\cdot D^\sigma\psi_1\,dx
	\right|.
	\end{split}
	\end{equation}
	The terms on the right-hand side can be 
	estimated using the Cauchy-Schwarz inequality and \eqref{pd-psi}
	(recall that $w=\partial_x u$, see Remark \ref{u_x=w}).
	Consequently, we successfully establish \eqref{energy-psi} 
	for $\psi_1$. The proof for $\psi_2$ proceeds similarly.
\end{proof}

\begin{remark}\label{C^1}
	From \eqref{ODE-phi-psi} and Proposition \ref{energy-phi-psi} 
	it follows that $(U,V)\in C^1(I_{\delta_0}, H^{s-2})$.
	Indeed, consider $(u,w,v,z)=(u,w,v,z)(x,t_0)$ and 
	$(\tilde{u},\tilde{w},\tilde{v},\tilde{z})
	=(u,w,v,z)(x,t)$ in \eqref{phi-psi} with $t,t_0\in I_{\delta_0}$.
	Taking into account \eqref{two-deriv-z} and \eqref{two-deriv-psi} 
	with $\sigma=s-2$, we conclude from \eqref{ODE-phi-psi} that 
	$\|\partial_t\phi_j\|_{H^{s-2}}$ and $\|\partial_t\psi_j\|_{H^{s-2}}$ 
	converge to zero as $t\to t_0$ for all $t_0\in I_{\delta_0}$ and $j=1,2$.
\end{remark}

Given Proposition \ref{energy-phi-psi}, we conclude that
\begin{equation*}
\left|
\frac{d}{dt}
\sum\limits_{i=1}^2\bigl(
\|\phi_i(t)\|_{H^{\sigma}}
+\|\psi_i(t)\|_{H^{\sigma}}\bigr)
\right|
\lesssim_{\sigma,s,\delta_0}
\left(\|u_0\|^2_{H^s}+\|v_0\|^2_{H^s}\right)
\sum\limits_{i=1}^2\bigl(
\|\phi_i(t)\|_{H^{\sigma}}
+\|\psi_i(t)\|_{H^{\sigma}}\bigr),
\quad t\in I_{\delta_0},
\end{equation*}
and therefore 
\begin{equation}\label{phi-psi-size}
\sum\limits_{i=1}^2\bigl(
\|\phi_i(t)\|_{H^{\sigma}}
+\|\psi_i(t)\|_{H^{\sigma}}\bigr)
\leq
\sum\limits_{i=1}^2\bigl(
\|\phi_i(0)\|_{H^{\sigma}}
+\|\psi_i(0)\|_{H^{\sigma}}\bigr)
e^{C|t|},\quad C>0,\,\,t\in I_{\delta_0}.
\end{equation}
Since $\phi_j(0)=\psi_j(0)=0$, $j=1,2$, the latter 
implies that $\phi_j(t)=\psi_j(t)=0$, $j=1,2$, for all $t\in I_{\delta_0}$.
Thus combining the results of this subsection and Subsection \ref{ex-sol-ODE}, 
we prove the existence and uniqueness of the solution of \eqref{ODE} in 
$C(I_{\delta_0}, H^{s-1})\cap C^1(I_{\delta_0}, H^{s-2})$.
Since $\partial_xu=w$ and $\partial_xv=z$ (see Remark \ref{u_x=w}), 
we conclude that $u,v\in C(I_{\delta_0}, H^{s})\cap C^1(I_{\delta_0}, H^{s-1})$.
Finally, since $\|u\|_{H^s}+\|v\|_{H^s}\leq\|(U,V)\|_{H^{s-1}}$,
\eqref{size-est-or} implies \eqref{size-est-u-v}. 
This, together with \eqref{ODE-a}--\eqref{ODE-d}, 
implies \eqref{deriv-est-u-v} and 
the proof of Theorem \ref{ex-un-sol} is complete.

\section{Continuity of the data-to-solution map}\label{Cont}
In this subsection we establish the continuous dependence 
of the solution of the two-component system \eqref{two-comp} 
on the initial data in $H^{s}$, $s>\frac{5}{2}$.

Let $u,v\in C(I_{\delta_0}, H^{s})\cap C^1(I_{\delta_0}, H^{s-1})$, 
$s>\frac{5}{2}$ be a local solution of \eqref{two-comp} 
with initial data $u(x,0)=u_0(x)$ and $v(x,0)=v_0(x)$.
Suppose that $u_{0,n}(x)\to u_0(x)$ and $v_{0,n}(x)\to v_0(x)$ in $H^{s}$ 
as $n\to\infty$. Then the goal is to prove that
\begin{subequations}
	\begin{align}
	\label{cont-dm-a}
	&\lim\limits_{n\to\infty}
	\left(\|u-u_n\|_{C(I_{\delta_0}, H^{s})}
	+\|v-v_n\|_{C(I_{\delta_0}, H^{s})}\right)=0,\\
	\label{cont-dm-b}
	&\lim\limits_{n\to\infty}
	\left(\|\partial_tu-\partial_tu_n\|_{C(I_{\delta_0}, H^{s-1})}
	+\|\partial_tv-\partial_tv_n\|_{C(I_{\delta_0}, H^{s-1})}\right)=0,
	\end{align}
\end{subequations}
which is equivalent to \eqref{cont-dep-lim} in Theorem \ref{cont-dep}.
Here the pair $u_n,v_n\in C(I_{\delta_0}, H^{s})\cap C^1(I_{\delta_0}, H^{s-1})$
is the solution of \eqref{two-comp} with initial data 
$u_{0,n},v_{0,n}$, specifying $\delta_0$ such that the
solutions $u_n,v_n$ and $u,v$ have the same lifespan $I_{\delta_0}$.

As in Section \ref{loc-ex-un}, we reduce two-component 
system \eqref{two-comp} to the system of 
PDEs \eqref{ODE-a}--\eqref{ODE-d} in $H^{s-1}$.
Then  
$(U,V)=(u,w,v,z)$ 
and
$(U_n,V_n)=(u_n,w_n,v_n,z_n)$, with 
$w_n=\partial_x u_n$, $z_n=\partial_x v_n$,
are the solutions of the Cauchy problem \eqref{ODE} with initial data 
$(u_0,w_0,v_0,z_0)$
and
$(u_{0,n},w_{0,n},v_{0,n},z_{0,n})$,
with
$w_{0,n}=\partial_x u_{0,n}$,
$z_{0,n}=\partial_x v_{0,n}$,
respectively.
Introduce two auxiliary solutions $(U^\ve,V^\ve)$ and 
$(U^\ve_n,V^\ve_n)$, $\ve\in(0,1)$ of \eqref{ODE} with mollified initial data
$(u_0^\ve,w_0^\ve,v_0^\ve,z_0^\ve)$
and
$(u_{0,n}^\ve,w_{0,n}^\ve,v_{0,n}^\ve,z_{0,n}^\ve)$, respectively.
Then, by the triangle inequality,
\begin{equation}\label{tr-cont}
\begin{split}
&\|(U,V)-(U_n,V_n)\|_{C(I_{\delta_0}, H^{s-1})}
\leq
\|(U,V)-(U^\ve,V^\ve)\|_{C(I_{\delta_0}, H^{s-1})}
\\ & \qquad +\|(U^\ve,V^\ve)-(U^\ve_n,V^\ve_n)\|_{C(I_{\delta_0}, H^{s-1})}
+\|(U^\ve_n,V^\ve_n)-(U_n,V_n)\|_{C(I_{\delta_0}, H^{s-1})}.
\end{split}
\end{equation}
To estimate the three terms on the right-hand side of \eqref{tr-cont}, 
we establish energy estimates 
for
$\|(U,V)-(U^\ve,V^\ve)\|_{H^{s-1}}$,
$\|(U^\ve_n,V^\ve_n)-(U_n,V_n)\|_{H^{s-1}}$
and
$\|(U^\ve,V^\ve)-(U^\ve_n,V^\ve_n)\|_{H^{s-1}}$.

The energy estimate for 
$\|(U^\ve,V^\ve)-(U^\ve_n,V^\ve_n)\|_{H^{s-1}}$
is the next lemma.

\begin{lemma}\label{en-est-ve-n}
	For all $t\in I_{\delta_0}$,
	\begin{equation}\label{en-est-2}
	\left|
	\frac{d}{dt}
	\|(U^\ve,V^\ve)(t)-(U^\ve_n,V^\ve_n)(t)\|_{H^{s-1}}
	\! \right|
	\lesssim_{s,\delta_0}\!
	\frac{1}{\ve}\bigl(\|u_0\|^2_{H^s}+\|v_0\|^2_{H^s}\bigr)
	\|(U^\ve,V^\ve)(t)
	-(U^\ve_n,V^\ve_n)(t)\|_{H^{s-1}}.
	\end{equation}
\end{lemma}

\begin{proof}
	Denote by $\phi_j$ and $\psi_j$, $j=1,2$, the 
	following differences (cf.\,\,\eqref{phi-psi}):
	\begin{equation}
	\nonumber
	\phi_1=u^\ve-u^\ve_n,\,\,
	\psi_1=w^\ve-w^\ve_n,\,\,
	\phi_2=v^\ve-v^\ve_n,\,\,
	\psi_2=z^\ve-z^\ve_n.
	\end{equation}
	The functions $\phi_j$, $\psi_j$ satisfy system \eqref{ODE-phi-psi} 
	with $(u^\ve,w^\ve,v^\ve,z^\ve)$
	and $(u^\ve_n,w^\ve_n,v^\ve_n,z^\ve_n)$ instead of 
	$(u,w,v,z)$ and $(\tilde{u},\tilde{w},\tilde{v},\tilde{z})$, respectively.
	To obtain the energy estimates for $\|\phi_j\|_{H^{s-1}}$ 
	($\|\psi_j\|_{H^{s-1}}$), we apply the operator $D^{s-1}$ 
	to the corresponding equation and multiply by 
	$D^{s-1}\phi_j$ ($D^{s-1}\psi_j$).
	All the integrals can be estimated using the 
	Kato-Ponce inequality, the algebra property, 
	the Sobolev theorem, and the Cauchy-Schwarz inequality, similarly 
	to the proof of Proposition \ref{energy-est}, except for those integrals 
	involving $\partial_x w_n^\ve$ and $\partial_x z_n^\ve$, 
	see $\beta_j$, $j=1,3,4,5$, and $\delta_j$, $j=1,2,3,5$.
	To handle these terms we use the inequality
	\begin{equation}\label{pd-moll}
	\|\partial_x w_n^\ve\|_{H^{s-1}}
	\lesssim\frac{1}{\ve}\|w_n\|_{H^{s-1}},
	\end{equation}
	which allows us to obtain the following estimate 
	(see the term $(\partial_x w_n^\ve)z_n^\ve$ in $\beta_1$):
	\begin{equation}
	\nonumber
	\begin{split}
	\left|
	\int_A D^{s-1}(\partial_x w_n^\ve)z_n^\ve\phi_1\cdot
	D^{s-1}\psi_1\,dx
	\right|
	&\leq\|(\partial_x w_n^\ve)z_n^\ve\phi_1\|_{H^{s-1}}
	\|\psi_1\|_{H^{s-1}}\\
	&\lesssim_s
	\frac{1}{\ve}\|w_n\|_{H^{s-1}}\|z_n^\ve\|_{H^{s-1}}
	\|\phi_1\|_{H^{s-1}}\|\psi_1\|_{H^{s-1}}.
	\end{split}
	\end{equation}
	All the other integral terms involving $\partial_x w_n^\ve$ and 
	$\partial_x z_n^\ve$ can be treated similarly.
	Finally, using the size estimate \eqref{size-est-or} 
	we arrive at \eqref{en-est-2}.
\end{proof}

Next, we derive the energy estimates for the two 
remaining terms on the right-hand side of \eqref{tr-cont}.
\begin{lemma}
	For $t\in I_{\delta_0}$,
	\begin{subequations}\label{en-est-13}
		\begin{align}
		\label{en-est-1}
		&\left|
		\frac{d}{dt}\|(U,V)(t)-(U^\ve,V^\ve)(t)\|_{H^{s-1}}
		\right|
		\\ & \nonumber 
		\qquad\qquad \lesssim_{s,\delta_0}
		\bigl(\|u_0\|^2_{H^s}+\|v_0\|^2_{H^s}\bigr)
		\|(U,V)(t)-(U^\ve,V^\ve)(t)\|_{H^{s-1}}+o(1),\\
		\label{en-est-3}
		&\left|
		\frac{d}{dt}\|(U^\ve_n,V^\ve_n)(t)
		-(U_n,V_n)(t)\|_{H^{s-1}}
		\right|
		\\ &\nonumber
		\qquad\qquad \lesssim_{s,\delta_0}
		\bigl(\|u_0\|^2_{H^s}+\|v_0\|^2_{H^s}\bigr)
		\|(U^\ve_n,V^\ve_n)(t)
		-(U_n,V_n)(t)\|_{H^{s-1}}+o(1),
		\end{align}
	\end{subequations}
	where $o(1)\to 0$ as $\ve\to0$.
\end{lemma}

\begin{proof}
	We give a proof for \eqref{en-est-1}, the 
	estimate \eqref{en-est-3} can be obtained similarly.
	Denote by $\phi_j$ and $\psi_j$, $j=1,2$ the following 
	differences (cf.\,\,\eqref{phi-psi})
	\begin{equation}
	\nonumber
	\phi_1=u-u^\ve,\,\,
	\psi_1=w-w^\ve,\,\,
	\phi_2=v-v^\ve,\,\,
	\psi_2=z-z^\ve.
	\end{equation}
	Following the methodology of Lemma \ref{en-est-ve-n}, we 
	apply the operator $D^{s-1}$ and multiply it by 
	either $D^{s-1}\phi_j$ or $D^{s-1}\psi_j$ in the corresponding equation. 
	In a similar manner as detailed in Proposition \ref{energy-est}, 
	we can estimate all integrals, except for those 
	involving $\partial_x w_n^\ve$ and $\partial_x z_n^\ve$, as seen 
	in $\beta_j$ for $j=1,3,4,5$ and $\delta_j$ for $j=1,2,3,5$. 
	To estimate these latter terms, note that (see the 
	term $u\partial_x w^\ve$ in $\beta_5$)
	\begin{equation}
	\nonumber
	\begin{split}
	\left|
	\int_A D^{s-1}u(\partial_x w^\ve)\psi_2
	\cdot D^{s-1}\psi_1\,dx
	\right|
	\leq&
	\left|
	\int_A [D^{s-1},u\psi_2]\partial_x w^\ve
	\cdot D^{s-1}\psi_1\,dx
	\right|\\
	&+
	\left|
	\int_A u\psi_2D^{s-1}\partial_x w^\ve
	\cdot D^{s-1}\psi_1\,dx
	\right|=I_{10}+I_{11}.
	\end{split}
	\end{equation}
	We estimate $I_{10}$ by the Kato-Ponce inequality, 
	while for $I_{11}$ we have
	\begin{equation}
	\nonumber
	\begin{split}
	I_{11}\leq\|u\psi_2\|_{L^\infty}
	\|\partial_xw^\ve\|_{H^{s-1}}
	\|\psi_1\|_{H^{s-1}}
	&\lesssim_s\frac{1}{\ve}\|u\|_{H^{s-2}}
	\|\psi_2\|_{H^{s-2}}
	\|w\|_{H^{s-1}}
	\|\psi_1\|_{H^{s-1}},\\
	&\lesssim_s o(1)\cdot\|u\|_{H^{s-1}}	\|w\|_{H^{s-1}}
	\|z\|_{H^{s-1}}\|\psi_1\|_{H^{s-1}},
	\end{split}
	\end{equation}
	where we have used \eqref{pd-moll} and the fact that
	$$
	\|\psi_2\|_{H^{s-2}}=\|(J_\ve-I)z\|_{H^{s-2}}
	\leq\|J_\ve-I\|_{\mathcal{L}(H^{s-1},H^{s-2})}
	\|z\|_{H^{s-2}}
	=o(\ve)\cdot\|z\|_{H^{s-2}}.
	$$
	All the other integrals involving $\partial_x w_n^\ve$ and 
	$\partial_x z_n^\ve$ can be treated similarly.
	Finally, using the size estimate \eqref{size-est-or}, we 
	arrive at \eqref{en-est-1}.
\end{proof}
By applying the Gronwall lemma to \eqref{en-est-2} 
and \eqref{en-est-13}, and considering the constraint 
$|t|\leq T_{\delta_0}$, we can deduce the subsequent 
inequalities valid for all $t$ in the interval $I_{\delta_0}$:
\begin{subequations}
	\begin{align}
	\label{appr-1}
	&\|(U^\ve,V^\ve)(t)-(U^\ve_n,V^\ve_n)(t)\|_{H^{s-1}}
	\leq\|(U^\ve,V^\ve)(0)-(U^\ve_n,V^\ve_n)(0)\|_{H^{s-1}}
	e^{\frac{C_{s,\delta_0}}{\ve}T_{\delta_0}},\\
	\label{appr-2}
	&\|(U,V)(t)-(U^\ve,V^\ve)(t)\|_{H^{s-1}}
	\leq\|(U,V)(0)-(U^\ve,V^\ve)(0)\|_{H^{s-1}}
	e^{C_{s,\delta_0}T_{\delta_0}}+o(1),\\
	\label{appr-3}
	&\|(U^\ve_n,V^\ve_n)(t)-(U_n,V_n)(t)\|_{H^{s-1}}
	\leq\|(U^\ve_n,V^\ve_n)(0)-(U_n,V_n)(0)\|_{H^{s-1}}
	e^{C_{s,\delta_0}T_{\delta_0}}+o(1).
	\end{align}
\end{subequations}
For any given $\hat\ve>0$, we can find a sufficiently 
small $\ve>0$ such that the right-hand sides of 
both \eqref{appr-2} and \eqref{appr-3} are less than $\hat\ve/3$ for all $n$. 
Consequently, for such a value of $\ve$, there exists a certain $N>0$ 
such that the right-hand side of equation \eqref{appr-1} is less 
than $\hat\ve/3$ for all $n\geq N$. Hence, using \eqref{tr-cont}, we 
can infer that $(U,V)$ converges to $(U_n,V_n)$ 
in ${C(I_{\delta_0}, H^{s-1})}$. This result, in turn, 
implies \eqref{cont-dm-a}.

To validate \eqref{cont-dm-b}, we take into account the 
system \eqref{ODE-phi-psi} with $(\tilde{u},\tilde{w},\tilde{v},\tilde{z})
=(u_n,w_n,v_n,z_n)$. Leveraging \eqref{two-deriv-z} 
and \eqref{two-deriv-psi} (see Remark \ref{C^1}), we can deduce 
that $\partial_t(U,V)$ converges to $\partial_t(U_n,V_n)$ 
in ${C(I_{\delta_0}, H^{s-2})}$. Hence, this completes the 
proof of Theorem \ref{cont-dep}.

\section{H\"older continuity properties of the solution}\label{Holder}
In this section, we explore the H{\"o}lder continuity 
properties of the data-to-solution map for the 
two-component system \eqref{two-comp} 
within the $H^r$ space for $r<s$. Initially, let us establish 
the subsequent lemma, which plays a crucial role in 
deriving the estimates detailed in Subsection \ref{Lip}.

\begin{lemma}\label{lemma-alg-r-s}
	Consider $(s,r)\in\mathbb{R}^2$ such that $s>\frac{5}{2}$, 
	$r\in[0,\frac{3}{2}]$, $r\geq3-s$ and $r<s-\frac{3}{2}$. 
	Then we have the following estimate (recall 
	that $A=\mathbb{R}$ or $A=\mathbb{R}/2\pi\mathbb{Z}$):
	\begin{equation}\label{alg-r-s}
	\|fg\|_{H^{r-1}(A)}\lesssim_{r,s}
	\|f\|_{H^{s-2}(A)}\|g\|_{H^{r-1}(A)}.
	\end{equation}
\end{lemma}

\begin{proof}
	First, consider $r\in[0,1]$.
	Since $s-1>\frac{3}{2}$ and 
	$r+(s-1)\geq 2$, we can apply \eqref{HH}
	with $(\hat{s},\hat{r})=(s-1,r)$,
	and thus arrive at \eqref{alg-r-s} for $r\in[0,1]$.
	
	Now let us prove \eqref{alg-r-s} for $r\in(1,\frac{3}{2}]$, 
	$r<s-\frac{3}{2}$ in the case $A=\mathbb{R}$ (cf.\,\,\cite[Lemma 2]{HH13}).
	Taking into account that 
	$\mathcal{F}(fg)=\mathcal{F}(f)*\mathcal{F}(g)$ we have
	\begin{equation*}
	\begin{split}
	\|fg\|^2_{H^{r-1}}&=
	\int_\mathbb{R}(1+k^2)^{r-1}\left|
	\int_\mathbb{R}\hat{f}(p)\hat{g}(k-p)\,dp
	\right|^2dk\\
	&=\int_\mathbb{R}(1+k^2)^{r-1}\left|
	\int_\mathbb{R}
	(1+p^2)^{\frac{r-\alpha}{2}}\hat{f}(p)
	(1+p^2)^{-\frac{r-\alpha}{2}}\hat{g}(k-p)\,dp
	\right|^2dk.
	\end{split}
	\end{equation*}
	By choosing $\alpha\geq 0$, employing the Cauchy-Schwarz 
	inequality for the inner integral, and switching the 
	order of integration, we deduce the following:
	\begin{equation*}
	\begin{split}
	\|fg\|^2_{H^{r-1}}&\leq
	\|f\|^2_{H^{r-\alpha}}
	\int_\mathbb{R}(1+k^2)^{r-1}
	\int_\mathbb{R}(1+p^2)^{\alpha-r}
	|\hat{g}(k-p)|^2\,dp\,dk\\
	&\leq\|f\|^2_{H^{r-\alpha}}	
	\int_\mathbb{R}(1+p^2)^{\alpha-r}
	\int_\mathbb{R}(1+k^2)^{r-1}
	|\hat{g}(k-p)|^2\,dk\,dp.
	\end{split}
	\end{equation*}
	Changing the variable $k\mapsto k-p$, using that 
	$1+(k+p)^2\leq2(1+k^2)(1+p^2)$ and taking 
	$\alpha\in[0,\frac{1}{2})$, we obtain
	\begin{equation*}
	\begin{split}
	\|fg\|^2_{H^{r-1}}&\lesssim
	\|f\|^2_{H^{r-\alpha}}
	\int_\mathbb{R}(1+p^2)^{\alpha-1}
	\int_\mathbb{R}(1+k^2)^{r-1}
	|\hat{g}(k)|^2\,dk\,dp\\
	&\lesssim\|f\|^2_{H^{r-\alpha}}	
	\|g\|^2_{H^{r-1}}
	\int_\mathbb{R}(1+p^2)^{\alpha-1}\,dp
	\lesssim_{\alpha}
	\|f\|^2_{H^{r-\alpha}}	
	\|g\|^2_{H^{r-1}}.
	\end{split}
	\end{equation*}
	For $s\in(\frac{5}{2},\frac{7}{2}]$ we use that $r<s-\frac{3}{2}$ 
	and take $\alpha=r-(s-2)$, while for 
	$s>\frac{7}{2}$ we use the estimate $\|f\|_{H^{r-\alpha}}\leq
	\|f\|_{H^{r}}\leq\|f\|_{H^{s-2}}$, and 
	thus \eqref{alg-r-s} is proved for $A=\mathbb{R}$.
	
	In the periodic case $A=\mathbb{R}/2\pi\mathbb{Z}$, 
	we have 
	$$
	f(x)=\sum\limits_{k\in\mathbb{Z}}
	e^{\mathrm{i}kx}\hat{f}(k), 
	\quad 
	\|f\|_{H^s}^2=2\pi\sum\limits_{k\in\mathbb{Z}}
	(1+k^2)^s|\hat{f}(k)|^2.
	$$
	Therefore, the demonstration of equation \eqref{alg-r-s} 
	for the periodic case remains identical, save for one distinction: 
	we replace the integral along the axis with a 
	summation over $k\in\mathbb{Z}$.
\end{proof}

\subsection{Lipschitz continuity for $(s,r)\in A_1$ and $(s,p)\in B_1$}\label{Lip}

Consider two solutions $u_j(x,t)$, $v_j(x,t)$, $j=1,2$, 
from Theorem \ref{Holder-cont}.
Then, as in Section \ref{loc-ex-un}, we reduce the Cauchy 
problems for the two-component system \eqref{two-comp} to the 
system \eqref{ODE-a}--\eqref{ODE-d}. 
Then subtracting the differential equations for 
$(u_1,w_1,v_1,z_1)$ and $(u_2,w_2,v_2,z_2)$, 
where $w_j=\partial_x u_j$ and $z_j=\partial_x v_j$, we 
obtain \eqref{ODE-phi-psi} with 
$(u,w,v,z)=(u_1,w_1,v_1,z_1)$ and
$(\tilde{u},\tilde{w},\tilde{v},\tilde{z})=(u_2,w_2,v_2,z_2)$.
Let us prove the following energy estimates for $\phi_j$ and $\psi_j$:
\begin{proposition}
	Set $I_{\rho,\delta_0}=
	[-T_{\rho,\delta_0},T_{\rho,\delta_0}]$. 
	Assume that $(s,r)\in A_1$. Then
	\begin{equation}
	\label{energy-psi-psi-H}
	\begin{split}
	&\left|\frac{d}{dt}\|\phi_j(t)\|_{H^{r-1}}
	\right|
	\lesssim_{r,s,\rho,\delta_0}
	\sum\limits_{i=1}^2\bigl(
	\|\phi_i(t)\|_{H^{r-1}}
	+\|\psi_i(t)\|_{H^{r-1}}\bigr),
	\quad t\in I_{\rho,\delta_0},\,\,j=1,2,\\
	&\left|\frac{d}{dt}\|\psi_j(t)\|_{H^{r-1}}
	\right|
	\lesssim_{r,s,\rho,\delta_0}
	\sum\limits_{i=1}^2\bigl(
	\|\phi_i(t)\|_{H^{r-1}}
	+\|\psi_i(t)\|_{H^{r-1}}\bigr),
	\quad t\in I_{\rho,\delta_0},\,\,j=1,2.
	\end{split}
	\end{equation}
\end{proposition}

\begin{proof}
	Let us first examine $(s,r)$ within the set $\{(s,r):s>5/2,r\in(3/2,s-1]\}$. 
	For such values of $r$, we can carry out the proof in a 
	manner similar to the methodology used in 
	Proposition \ref{energy-phi-psi}, with $\sigma$ 
	equated to $r-1$.

	Next, let us turn our attention to $(s,r)$ in the set 
	$\{(s,r):s>5/2,r\in[0,3/2],r\in[3-s,s-3/2)\}$. Here, for this 
	range of $r$, it is noted that the algebra property fails 
	in $H^{r-1}$. However, by applying arguments akin to 
	those employed in Proposition \ref{energy-phi-psi} 
	with $\sigma=r-1$ and substituting the algebra property 
	estimate with $\eqref{alg-r-s}$, we can still successfully 
	arrive at \eqref{energy-psi-psi-H}.
\end{proof}

From \eqref{energy-psi-psi-H} it follows that (cf.\,\,\eqref{phi-psi-size})
\begin{equation}\label{phi-psi-size-H}
\sum\limits_{i=1}^2\bigl(
\|\phi_i(t)\|_{H^{r-1}}
+\|\psi_i(t)\|_{H^{r-1}}\bigr)
\lesssim_{r,s,\rho,\delta_0}
\sum\limits_{i=1}^2\bigl(
\|\phi_i(0)\|_{H^{r-1}}
+\|\psi_i(0)\|_{H^{r-1}}\bigr),\,\,t\in I_{\rho,\delta_0}.
\end{equation}
Given that $\psi_j=\partial_x\phi_j$, we can derive 
\eqref{H-cont-sol} from \eqref{phi-psi-size-H} for all 
pairs $(s,r)$ belonging to the set $A_1$.

\begin{remark}\label{Lip-A-1}
	In this study, we establish the Lipschitz property 
	of the data-to-solution map for the two-component 
	system \eqref{two-comp}. Our analysis covers a more 
	expansive subset of the $(s,r)$-plane than was 
	previously demonstrated in \cite{HM14-1} for the FORQ equation. 
	This can be seen in Figure \ref{A-regions} and \cite[Figure 1]{HM14-1}, 
	specifically regarding the subset $A_1$ in \eqref{A-j-regions} 
	and the definition of $A_1$ in \cite[Theorem 1]{HM14-1}. 
	Our proof utilizes inequality \eqref{alg-r-s} for $r\in(1/2,3/2]$, which 
	provides a more precise estimation than the 
	inequality $\|fg\|_{H^r}\lesssim_{r} \|f\|_{H^{r+1}}\|g\|_{H^r}$ 
	used in \cite{HM14-1}.
\end{remark}

Let us prove the Lipschitz property for $(s,p)\in B_1$.
Consider a pair $(s,\hat{r})\in\tilde{A}_1$, where
$\tilde{A}_1$ has the form
$$
\tilde{A}_1=\{(s,\hat{r}):s>5/2,\,
\hat r\in(5/2,s]\}\cup
\{(s,\hat r):s>5/2,\,
\hat r\in[1,5/2],\,\hat r\in[4-s,s-1/2)\}.
$$
Next, we utilize either the algebra property or 
Lemma \ref{lemma-alg-r-s}, contingent on the value of $\hat r$. 
From \eqref{ODE-phi-psi}, we then derive the 
subsequent estimate:
\begin{equation}\label{pd-t-phi-psi}
\sum\limits_{i=1}^2\bigl(
\|\partial_t\phi_i(t)\|_{H^{\hat r-2}}
+\|\partial_t\psi_i(t)\|_{H^{\hat r-2}}\bigr)
\lesssim_{\hat r,s,\rho}
\sum\limits_{i=1}^2\bigl(
\|\phi_i(t)\|_{H^{\hat r-1}}
+\|\psi_i(t)\|_{H^{\hat r-1}}\bigr).
\end{equation}
Given \eqref{phi-psi-size-H}, we can infer 
from \eqref{pd-t-phi-psi} that the following 
estimate holds true for all pairs $(s,\hat r)$ that are 
part of $A_1\cap \tilde{A}_1$:
$$
\sum\limits_{i=1}^2\bigl(
\|\partial_t\phi_i(t)\|_{H^{\hat r-2}}
+\|\partial_t\psi_i(t)\|_{H^{\hat r-2}}\bigr)
\lesssim_{\hat r,s,\rho,\delta_0}
\sum\limits_{i=1}^2\bigl(
\|\phi_i(0)\|_{H^{\hat r-1}}
+\|\psi_i(0)\|_{H^{\hat r-1}}\bigr),\,\,t\in I_{\rho,\delta_0}.
$$
Recalling that $\psi_j=\partial_x\phi_j$, $j=1,2$, and 
specifying $p=\hat r-1$, we arrive at \eqref{H-cont-pd-t} 
for $(s,p)\in B_1$.

\subsection{H\"older continuity for $(s,r)\in A_i$ and $(s,p)\in B_j$, $i,j>1$}

To prove the H\"older property for the rest of the 
regions we will use the interpolation inequality \eqref{interpolation}.

Let us prove the H\"older continuity in the regions $A_i$, $i=2,3,4,5$. 
For $(s,r)\in A_2$, we have that $r<3-s<s$ and using \eqref{interpolation} with 
$(\sigma_1,\sigma,\sigma_2)=(r,3-s,s)$ 
as well as the Lipschitz property for $r=3-s$, $s\in(5/2,3]$, we obtain 
(recall that $\|u_{0,j}\|_{H^s}, 
\|v_{0,j}\|_{H^s}\leq\rho$, $j=1,2$)
\begin{equation}\label{A-2}
\begin{split}
\|(u_1,v_1)-(u_2,v_2)\|_{H^r}
&\leq \|(u_1,v_1)-(u_2,v_2)\|_{H^{3-s}}
\lesssim_{r,s,\rho,\delta_0}
\|(u_{0,1},v_{0,1})
-(u_{0,2},v_{0,2})\|_{H^{3-s}}\\
&\lesssim_{r,s,\rho,\delta_0}
\|(u_{0,1},v_{0,1})
-(u_{0,2},v_{0,2})\|_{H^{r}}^{\frac{2s-3}{s-r}}.
\end{split}
\end{equation}

In the case $(s,r)\in A_4$ and $(s,r)\in A_6$ we 
have $r<3/2+\ve_0/2 <s$ and $r<0<s$, respectively. 
Using \eqref{interpolation} with 
$(\sigma_1,\sigma,\sigma_2)=(r,3/2+\ve_0/2,s)$
and
$(\sigma_1,\sigma,\sigma_2)=(r,0,s)$, respectively, and 
arguing as in \eqref{A-2}, we arrive at 
\eqref{H-cont-sol} for $(s,r)\in A_4\cup A_6$.

Consider $(s,r)\in A_5$. In this case, $s-1<r<s$ and thus, 
by \eqref{interpolation} with 
$(\sigma_1,\sigma,\sigma_2)=(s-1,r,s)$,
we have
$$
\|(u_1,v_1)-(u_2,v_2)\|_{H^r}
\leq
\|(u_1,v_1)-(u_2,v_2)\|_{H^{s-1}}^{s-r}
\|(u_1,v_1)-(u_2,v_2)\|_{H^s}^{s-r+1}.
$$
Using the size estimate \eqref{size-est-u-v} as well 
as Lipschitz property for $r=s-1$, $s>\frac{5}{2}$, 
from this we deduce 
$$
\|(u_1,v_1)-(u_2,v_2)\|_{H^r}
\lesssim_{\rho,\delta_0}
\|(u_1,v_1)-(u_2,v_2)\|_{H^{r}}^{s-r}
\lesssim_{r,s\rho,\delta_0}
\|(u_{0,1},v_{0,1})
-(u_{0,2},v_{0,2})\|_{H^{r}}^{s-r},
$$
which implies \eqref{H-cont-sol} for $(s,r)\in A_5$.

In conclusion, for the case where $(s,r)$ falls within 
the set $A_3$, we observe that $s-3/2-\ve_0/2<r<s$. 
By applying \eqref{interpolation} with the 
trio $(\sigma_1,\sigma,\sigma_2) = (s-3/2-\ve_0/2,r,s)$ and 
reasoning along the lines of our approach 
for the region $A_5$, we derive \eqref{H-cont-sol} 
for all $(s,r)\in A_3$.

\begin{remark}\label{Hold-exp-3-4}
	It is important to note that the H{\"o}lder exponent 
	$\gamma=\gamma(s,r)$, as observed in our study, presents 
	a more optimal outcome within the regions $A_3\cup A_4$ 
	compared to the results derived in \cite{HM14-1}. 
	This can be particularly contrasted with the H{\"o}lder 
	exponent $\alpha$, defined in \cite[(1.17)]{HM14-1} 
	for the region $A_3$, which is notably not 
	continuous along $r=3/2$. In our research, we find that within 
	$A_3\cup A_4$, we can apply \eqref{interpolation} under the 
	conditions of either $s-3/2-\ve_0/2<r<s$ or $r<3/2+\ve_0/2 <s$. 
	By evaluating the maximum value between $\frac{2(s-r)}{3+\ve_0}$ 
	and $\frac{2s-3-\ve_0}{2(s-r)}$ for all pairs $(s,r)$ within $A_3\cup A_4$, 
	we derive the H{\"o}lder exponent $\gamma$ for these regions.
\end{remark}

Now, we aim to establish the H{\"o}lder continuity 
within the regions $B_j$, for $j=2,\dots,6$. For all pairs $(s,p)$ that 
belong to the regions $B_3$, $B_4$, and $B_6$, we adopt 
a similar approach to the one used for regions $A_2$, $A_4$, and $A_6$. 
We assign $(\sigma_1,\sigma,\sigma_2)$ as follows:

\smallskip

for $B_3$, $(\sigma_1,\sigma,\sigma_2) = (p,1/2+\ve_1/2,s-1)$;

for $B_4$, $(\sigma_1,\sigma,\sigma_2) = (p,3-s,s-1)$;

and for $B_6$, $(\sigma_1,\sigma,\sigma_2) = (p,0,s-1)$.

\smallskip

\noindent By following this approach, we arrive 
at \eqref{H-cont-pd-t} for these regions.

Conversely, for the regions $B_2$ and $B_5$, where $(s,p)$ 
belongs, we employ the estimation methods used for 
regions $A_3$ and $A_5$. We set $(\sigma_1,\sigma,\sigma_2)$ to:

\smallskip

for $B_2$, $(\sigma_1,\sigma,\sigma_2) = (s-5/2-\ve_1/2,p,s-1)$;

and for $B_5$, $(\sigma_1,\sigma,\sigma_2) = (s-2,p,s-1)$.

\smallskip

\noindent By adhering to this process, we 
ultimately derive \eqref{H-cont-pd-t}.


\begin{thebibliography}{99}
	\bibitem{AKNS74}
	M.J. Ablowitz, D.J. Kaup, A.C. Newell and H. Segur.
	\newblock The inverse scattering transform-Fourier analysis 
	for nonlinear problems.
	\newblock {\em Stud. Appl. Math.}, 53:249--315, 1974.
	
	\bibitem{AM13}
	M.J. Ablowitz and Z.H. Musslimani.
	\newblock Integrable nonlocal nonlinear 
	Schr{\"o}dinger equation.
	\newblock {\em Phys. Rev. Lett.}, 110:064105, 2013.
	
	\bibitem{AM17}
	M.J. Ablowitz and Z.H. Musslimani.
	\newblock Integrable nonlocal nonlinear equations.
	\newblock {\em Stud. Appl. Math.}, 139:7--59, 2017.
	
	\bibitem{BB98}
	C.M. Bender, S. Boettcher.
	\newblock Real spectra in non-Hermitian Hamiltonians having P-T symmetry.
	\newblock {\em Phys. Rev. Lett.}, 80:5243, 1998.
	
	\bibitem{BBCF07}
	C.M. Bender, D.C. Brody, J.-H. Chen and E. Furlan.
	\newblock 
	\textit{PT}-symmetric 
	extension of the Korteweg-de Vries equation.
	\newblock {\em J. Phys. A: Math. Theor.}, 40:F153, 2007.
	
	\bibitem{BKS20}
	A. Boutet de Monvel, I. Karpenko and D. Shepelsky.
	\newblock A {R}iemann-{H}ilbert approach to the modified
	Camassa-Holm equation with nonzero boundary
	conditions.
	\newblock{\em J. Math. Phys.}, 61:031504, 2020.
	
	\bibitem{BC07}
	A. Bressan and A. Constantin.
	\newblock Global conservative solutions of the Camassa-Holm equation.
	\newblock {\em Arch. Ration. Mech. Anal.}, 183:215--239, 2007.
	
	\bibitem{D01}
	R. Danchin.
	\newblock A few remarks on the Camassa-Holm equation.
	\newblock {\em Diff. Integral Eqns.}, 14:953--988, 2001.
	
	\bibitem{D69}
	J. Dieudonn\'{e}.
	\newblock Foundations of Modern Analysis.
	\newblock Academic Press, New York, 1969.
	
	\bibitem{EMK18}
	R. El-Ganainy, K. Makris, M. Khajavikhan, Z. Musslimani, 
	S. Rotter and D. Christodoulides.
	\newblock Non-Hermitian physics and PT symmetry.
	\newblock {\em Nature Phys.}, 14:11--19, 2018.
	
	\bibitem{F95}
	A.S. Fokas.
	\newblock On a class of physically important integrable equations.
	\newblock {\em Physica D}, 87:145--150, 1995.
	
	\bibitem{FGLQ13}
	Y. Fu, G. Gui, Y. Liu, and C. Qu.
	\newblock On the Cauchy problem for the integrable 
	modified Camassa-Holm equation with cubic nonlinearity.
	\newblock {\em J. Differ. Equations}, 255(7):1905--1938 (2013).
	
	\bibitem{FF81}
	B. Fuchssteiner and A.S. Fokas.
	\newblock Symplectic structures, their B\"acklund 
	transformations and hereditaries.
	\newblock {\em Physica D}, 4:47--66, 1981.
	
	\bibitem{GLOQ13}
	G. Gui, Y. Liu, P.J. Olver and C. Qu.
	\newblock Wave-breaking and peakons for a modified
	Camassa-Holm equation.
	\newblock {\em Commun. Math. Phys.}, 319:731--759, 2013.
	
	\bibitem{GLL18}
	Y. Gao, L. Li and J.-G. Liu.
	\newblock A Dispersive Regularization for the Modified 
	Camassa-Holm Equation.
	\newblock {\em SIAM J. Math. Anal.}, 50(3):2807--2838, 2018.
	
	\bibitem{HH13}
	A. Himonas, J. Holmes.
	\newblock 
	H\"older continuity of the solution map for the Novikov equation.
	\newblock {\em J. Math. Phys.}, 54:061501, 2013.
	
	\bibitem{HM14}
	A. Himonas, D. Mantzavinos. 
	The Cauchy problem for the Fokas-Olver-Rosenau-Qiao equation.
	\newblock {\em J. Nonlinear Analysis: Theory, Methods 
	\& Applications}, 95:499--529, 2014.

	\bibitem{HM14-1}
	A. Himonas, D. Mantzavinos.
	\newblock 
	H\"older continuity for the Fokas-Olver-Rosenau-Qiao equation.
	\newblock {\em J. Nonlinear. Sci.}, 24:1105--1124, 2014.
	
	\bibitem{HK09}
	A. Himonas and C. Kenig.
	\newblock Non-uniform dependence on initial data for the CH equation on the line.
	\newblock {\em Diff. Integral Eqns.}, 22:201--24, 2009.
	
	\bibitem{HR07}
	H. Holden and X. Raynaud, 
	\newblock Global conservative solutions of the 
	Camassa-Holm equation -- a Lagrangian point of view.
	\newblock {\em Comm. Partial Differential Equations}, 32:1511--1549, 2007.
	
	\bibitem{HFQ17}
	Y. Hou, E. Fan, and Z. Qiao.
	\newblock The algebro-geometric solutions for the 
	Fokas-Olver-Rosenau-Qiao (FORQ) hierarchy. 
	\newblock {\em J. Geom. Phys.}, 117:105--133, 2017.
	
	\bibitem{KP88}
	T. Kato and G. Ponce.
	\newblock Commutator estimates and the euler and navier-stokes equations.
	\newblock {Comm. Pure Appl. Math.}, 41(7):891--907, 1988.
	
	\bibitem{KLOQ16}
	J. Kang, X. Liu, P. Olver and C. Qu.
	\newblock Liouville correspondence between the modified KdV
	hierarchy and its dual integrable hierarchy.
	\newblock {\em J. Nonlinear Sci.}, 26:141--170, 2016.
	
	\bibitem{KYZ16}
	V.V. Konotop, J. Yang and D.A. Zezyulin.
	\newblock Nonlinear waves in PT-symmetric systems.
	\newblock {\em Rev. Mod. Phys.}, 88:035002, 2016.
	
	\bibitem{LLQ14}
	X. Liu, Y. Liu and C. Qu.
	\newblock Orbital stability of the train of peakons for an 
	integrable modified Camassa-Holm equation.
	\newblock {\em Adv. Math.}, 255:1--37, 2014.
	
	\bibitem{LH17}
	S.Y. Lou, F. Huang.
	\newblock Alice-Bob physics{:} coherent solutions of nonlocal KdV systems.
	\newblock {\em Sci. Rep.}, 7:869, 2017.

	\bibitem{LQ17}
	S.Y. Lou, Z. Qiao.
	\newblock Alice-Bob peakon systems.
	\newblock {\em Chin. Phys. Lett.}, 34(10):100201, 2017.
	
	\bibitem{LM21}
	L. Ling, W.X. Ma.
	\newblock Inverse scattering and soliton solutions of 
	nonlocal complex reverse-spacetime modified Korteweg-de Vries hierarchies.
	\newblock {\em Symmetry}, 13:512 2021.
	
	\bibitem{M78}
	F. Magri.
	\newblock A simple model of the integrable Hamiltonian equation.
	\newblock {\em J. Math. Phys.}, 19:1156--1162, (1978).
	
	\bibitem{MM13}
	Y. Mi and C. Mu.
	\newblock Well-posedness and analyticity for an integrable 
	two-component system with cubic nonlinearity.
	\newblock {\em J. Hyperbolic Differ. Equations}, 10(04):703--723, 2013.
	
	\bibitem{OR96}
	P.J. Olver and P. Rosenau.
	\newblock Tri-Hamiltonian duality between solitons and 
	solitary-wave solutions having compact support.
	\newblock {\em Phys. Rev. E}, 53:1900--1906, 1996.
	
	\bibitem{Q06}
	Z.J. Qiao.
	\newblock A new integrable equation with cuspons and W/M-shape-peaks solitons.
	\newblock {\em J. Math. Phys.}, 47:112701, 2006.
	
	\bibitem{SW04}
	T. Sch\"afer, C.E.Wayne.
	\newblock Propagation of ultra-short optical pulses in cubic nonlinear media. 
	\newblock {\em Physica D}, 196:90--105, 2004.
	
	\bibitem{SQQ11}
	J.F. Song, C.Z. Qu and Z.J. Qiao.
	\newblock A new integrable two-component system with cubic nonlinearity.
	\newblock {\em J. Math. Phys.}, 52:013503, 2011.
	
	\bibitem{TLH18}
	X.-Y. Tang, Z.-F. Liang and X.-Z. Hao.
	\newblock Nonlinear waves of a nonlocal modified KdV equation 
	in the atmospheric and oceanic dynamical system.
	\newblock {\em Commun. Nonlinear Sci. Numer. Simul.}, 60:62--71, 2018.
	
	\bibitem{T91}
	M. Taylor.
	\newblock Pseudodifferential Operators and Nonlinear PDE.
	\newblock Boston, MA: Birkhauser.

	\bibitem{T03}
	M. Taylor.
	\newblock Commutator estimates.
	\newblock {\em Proc. Am. Math. Soc.}, 131:1501--1507, 2003.
	
	\bibitem{TL13}
	K. Tian and Q.P. Liu.
	\newblock Tri-Hamiltonian duality between the 
	Wadati-Konno-Ichikawa hierarchy and the Song-Qu-Qiao hierarchy.
	\newblock {\em J. Math. Phys.}, 54:043513 2013.
	
	\bibitem{T78}
	H. Triebel. 
	\newblock Interpolation Theory, Function Spaces, Differential Operators.
	\newblock Amsterdam, North-Holland Pub. Co., 1978.
	
	\bibitem{T89}
	G. Tu.
	\newblock The trace identity, a powerful tool for constructing 
	the Hamiltonian structure of integrable systems.
	\newblock {\em J. Math. Phys.}, 30:330--338, 1989.
	
	\bibitem{WKI79}
	M. Wadati, K. Konno, and Y.H. Ichikawa.
	\newblock New integrable nonlinear evolution equations. 
	\newblock {\em J. Phys. Soc. Jpn.}, 47:1698--1700, 1979.
	
	\bibitem{WY23}
	Z. Wang and K. Yan.
	\newblock Blow-up data for a two-component 
	Camassa-Holm system with high order nonlinearity.
	\newblock {\em J. Differential Equations}, 358(15):256--294, 2023.
	
	\bibitem{XQZ15}
	B. Xia, Z. Qiao, R. Zhou.
	\newblock A synthetical two-component model with peakon solutions.
	\newblock \textit{Stud. Appl. Math.}, 135(3):248--276, 2015.
	
	\bibitem{YQZ15}
	K. Yan, Z. Qiao, and Y. Zhang. 
	\newblock Blow-up phenomena for an integrable two-component 
	Camassa-Holm system with cubic nonlinearity and peakon solutions.
	\newblock {\em J. Differential Equations}, 259(11):6644--6671, 2015.

\end{thebibliography}
\end{document}